\documentclass[12pt,english]{article}

\usepackage{hyperref}
\hypersetup{
    unicode = false,
    pdftoolbar = true,
    pdfmenubar = true,
    pdffitwindow = true,
    pdftitle = {Gengrad},
    pdfauthor = {Audet Hare},
    pdfsubject = {Gengrad},
    pdfnewwindow = true,
    pdfkeywords = {Gengrad},
    colorlinks = true,
    linkcolor = blue,
    citecolor = blue,
    filecolor = black,
    urlcolor = blue,
    breaklinks = true
}

\usepackage{graphicx}
\usepackage[english]{babel}
\usepackage{amsmath}
\usepackage{amssymb}
\usepackage{multirow}
\usepackage{epsfig}
\usepackage{epic}
\usepackage{pgf}
\usepackage{array}
\usepackage{geometry}
\usepackage{times}
\usepackage{xspace}
\usepackage{tikz}
\usetikzlibrary{arrows}

\usepackage{caption}
\usepackage{subcaption}

\def\R{{\mathbb{R}}}

\def\H{{\mathcal{H}}}
\def\ZZ{o_c}

\newcommand{\nf}{n_f} 		
\newcommand{\nv}{n_v}
\newcommand{\e}{{\textmd{e}}}
\newcommand{\ej}{{\textmd{e}_j}}
\newcommand{\en}{{\textmd{e}_n}}
\newcommand{\proof}{\noindent {\bf Proof: }}
\newcommand{\argmin}{\mathop{\mathrm{argmin}}}

\newcommand{\conv}{\mathop{\mathrm{conv}}}

\newtheorem{Def}{Definition}[section]
\newtheorem{Lem}[Def]{Lemma}
\newtheorem{Rem}[Def]{Remark}
\newtheorem{Prop}[Def]{Proposition}
\newtheorem{Th}[Def]{Theorem}
\newtheorem{Cor}[Def]{Corollary}
\newtheorem{Ex}[Def]{Example}

\newcommand{\oracle}{{\sf D}}
\newcommand{\dist}{{\mathrm{dist}}}
\newcommand{\qed}{\hrulefill$\diamondsuit$ \medskip}

\newtheorem{algolocal}[Def]{Algorithm}
\newenvironment{algo}[1]
	{	\begin{algolocal} {\em \bf #1}  \end{algolocal}
		\begin{list}{}{} \item \vspace*{-3mm}			}
	{	\end{list} \vspace*{-6mm}\noindent \hrulefill \\			}

\title{
    Algorithmic construction of the subdifferential from directional derivatives
        \thanks{
	Work of the first author was supported by Natural Sciences and Engineering Research Council of Canada (NSERC), Discovery grant  \#2015-05311.
    Work of the second author was supported by NSERC, Discovery Grant \#355571-2013.
    }
}

\author{
    \href{mailto:Charles.Audet@gerad.ca}{Charles Audet}\thanks{{GERAD}
          and D\'epartement de math\'ematiques et g\'enie industriel,
          \'Ecole Polytechnique de Montr\'eal,
          C.P. 6079, Succ. Centre-ville,
          Montr\'eal, Qu\'ebec, Canada H3C 3A7,
          \href{https://www.gerad.ca/Charles.Audet}{www.gerad.ca/Charles.Audet}.
  }
  \and
        \href{mailto:warren.hare@ubc.ca}{Warren Hare}\thanks{Mathematics, Irving K. Barber School, University of British Columbia, Okanagan Campus, Kelowna, B.C., V1V 1V7, Canada, \href{http://math.ok.ubc.ca/faculty/hare.html}{math.ok.ubc.ca/faculty/hare}.}
}

\begin{document}
\maketitle

\noindent
\begin{abstract}
The subdifferential of a function is a generalization for nonsmooth functions of the concept of gradient. It is frequently used in variational analysis, particularly in the context of nonsmooth optimization.  The present work proposes algorithms to reconstruct a polyhedral subdifferential of a function from the computation of finitely many directional derivatives. We provide upper bounds on the required number of directional derivatives when the space is $\R^1$ and $\R^2$, as well as in  $\R^n$ where subdifferential is known to possess at most three vertices.
\end{abstract}

\noindent
{\bf Key words.} subdifferential, directional derivative, polyhedral construction, Geometric Probing. \\

\noindent
{\bf AMS subject classifications.} Primary 52B12, 65K15; Secondary 49M37, 90C30, 90C56.

%

\newpage
\section{Introduction}\label{introduction}

The {\em subdifferential} of a nonsmooth function represents the set of generalized gradients for the function (a formal definition appears in Section \ref{Definitions}).  It can be used to detect descent directions \cite[Thm 8.30]{rockafellar-wets-1998}, check first order optimality conditions \cite[Thm 10.1]{rockafellar-wets-1998}, and create cutting planes for convex functions \cite[Prop 8.12]{rockafellar-wets-1998}.  It shows up in numerous algorithms for nonsmooth optimization: steepest descent \cite[Sec XII.3]{Hiriart-Urruty-Lemarechal-1993b}, projective subgradient methods \cite[Sec XII.4]{Hiriart-Urruty-Lemarechal-1993b}, bundle methods \cite[Sec XIV.3]{Hiriart-Urruty-Lemarechal-1993b}, etc.  Its calculus properties have been well researched \cite[Chpt 10]{rockafellar-wets-1998}, and many favourable rules have been developed.  Overall, it is reasonable to state that the subdifferential is one of the most fundamental objects in variational analysis.

Given its role in nonsmooth optimization, it is no surprise that some researchers have turned their attention to the question of how to approximate a subdifferential using `simpler information'.  Besides the mathematical appeal of such a question, such research has strong links to the fields of {\em Derivative-free Optimization} and {\em Geometric Probing}.

Derivative-free optimization focuses on the development of algorithms to minimize a function using only function values.  Thus, in this case, `simpler information' takes the form of function values.  The ability to use function evaluations to approximate generalized gradient information is at the heart of the convergence analyses of many derivative-free optimization methods~\cite{AuDe03a,AuDe2006,FiKe2004a,FiKe09}.  Some research has explicitly proposed methods to approximate subdifferential for the purposes of derivative-free optimization methods \cite{bagirov-2003, bagirov-karasozen-sezer-2008, Kiwiel-2010, hare-macklem-2012, Hare-Nutini-2013}.  Many of these researchers focus on approximating the subdifferential by using a collection of approximate {\em directional derivatives} \cite{bagirov-2003, bagirov-karasozen-sezer-2008, Kiwiel-2010}.
Recall that directional derivative provides the slope for a function in a given direction in the classical limiting sense of single variable calculus -- a formal definition appears in Section \ref{Definitions}. Directional derivatives are intimately linked to subdifferential maps (see Section \ref{Definitions}) and are appealing in that they can be approximated by simple finite difference formulae.
This makes them an obvious tool to approximate subdifferentials.
The present paper studies how many directional derivatives are needed to reconstruct the subdifferential.

Geometric Probing considers problems of determining the geometric structure of a set by using a probe \cite{Skiena89, Romanik95}.  If the set is a subdifferential and the probe is a directional derivative, then the Geometric Probing problem is to reconstruct the subdifferential using the `simpler information' of directional derivatives (details appear in Subsection \ref{Geoprobing}).  Geometric Probing first arose in the area of robotics, where tactile sensing is used to determine the shape of an object \cite{ColeYap84}.   The problem has since been well studied in $\R^2$ \cite{ColeYap84, Lindenbaum92} and partly studied in $\R^3$ \cite{Imiya2012}.  As Geometric Probing principally arises in robotics and computer vision, it is not surprising that literature outside of $\R^2$ and $\R^3$ appears absent.

In this paper we examine the links between directional derivatives and the ability to use them to reconstruct a subdifferential.  We focus on the easiest case, where the subdifferential is a polytope and the directional derivatives are exact.
Let $\nv$ denote the number of vertices of the subdifferential, and $\nf$ be a given upper bound on $\nv$.  (The value $\nf = \infty$ is accepted to represent the situation when no upper bound is known.)
We show that,
\begin{itemize}
    \item in $\R^1$, the subdifferential can be reconstructed using a single
     directional derivative evaluation if $\nf=1$,
     and $2$ evaluations otherwise (Subsection \ref{OneDim});
    \item in $\R^2$, if $\nf= \nv$,
    then the subdifferential can be reconstructed using $3\nv$ directional derivative evaluations (Subsection \ref{TwoDim});
    \item in $\R^2$, if $\nf > \nv$, then the subdifferential can be reconstructed using $3\nv+1$ directional derivative evaluations (Subsection \ref{TwoDim});
    \item in $\R^n$, if  $\nf = 2$, then the subdifferential can be reconstructed using $3n-1$ directional derivative evaluations (Subsection \ref{nDim-nf2}); and
    \item in $\R^n$, if  $\nf =3$ vertices, then the subdifferential can be reconstructed using $5n-1$ directional derivative evaluations (Subsection \ref{nDim-nf3}).
\end{itemize}
These results can be loosely viewed as providing a lower bound on the number of approximate directional derivative evaluations that would be required to create a good approximation of the subdifferential in derivative-free optimization.  The results also advance research in Geometric Probing, which historically has only considered polytopes in $\R^2$ and $\R^3$.

Before proving these results, Section \ref{Definitions} proposes the necessary background for this work.
It also includes a method of problem abstraction which links the research to Geometric Probing.
Sections~\ref{OneDim} and~\ref{TwoDim} present the results in $\R^1$ and in $\R^2$, and Section~\ref{nDim} focuses on $\R^n$.
The paper concludes with some thoughts on the challenge of reconstructing polyhedral subdifferentials when directional derivatives are only available via finite difference approximations, and some other possible directions for future research.

\section{Definitions and problem abstraction}\label{Definitions}

Given a nonsmooth function $f : \R^n \mapsto \R\cup\{+\infty\}$ and a point $\bar{x} \in \R^n$ where $f(\bar{x})$ is finite, we define the {\em regular subdifferential} of $f$ at $\bar{x}$, ${\partial} f(\bar{x})$, as the set
    $${\partial} f(\bar{x}) = \left\{ v \in \R^n ~:~ f(x) \geq f(\bar{x}) + v^\top(x-\bar{x}) + o\|x-\bar{x}\|\right\}.$$
If $f$ is convex, then this is equivalent to the classical subdifferential of convex analysis \cite[Prop 8.12]{rockafellar-wets-1998}, i.e.,
    $$\mbox{if $f$ is convex, then~} \partial f(\bar{x}) = \{ v \in \R^n ~:~ f(x) \geq f(\bar{x}) + v^\top(x-\bar{x}) ~\mbox{ for all } x \in \R^n\}.$$
In this paper, we consider the situation where the subdifferential is a polyhedral set.  This arises, for example, when $f$ is a {\em finite max function}.  In particular, if
\begin{equation}
	\label{eq-finitemax}
	f = \max \left\{f_i ~:~ i \in \{  1, 2, ..., \nf\}  \right\},
\end{equation}
where each $f_i \in \mathcal{C}^1$, then
\begin{equation}\label{eq:partconv}
    \partial f(\bar{x}) = \conv\{ \nabla f_i(\bar{x}) ~:~ i \in A(\bar{x})\}
\end{equation}
where $A(\bar{x}) = \{ i \in \{  1, 2, ..., \nf\} : f_i(\bar{x}) = f(\bar{x})\}$ \cite[Ex 8.31]{rockafellar-wets-1998}. It follows that, in this case, the subdifferential $\partial f(\bar{x})$ is a (nonempty) polytope with at most $\nf$ vertices.

Related to the subdifferential is the {\em directional derivative}.
Formally,
	the directional derivative $df(\bar{x};\bar{d})$ of a continuous function $f : \R^n \mapsto \R$
	at $\bar x \in \R^n$ in the direction $\bar d \in \R^n$ is
$$df(\bar{x};\bar{d}) = \lim_{\tau \searrow 0} \frac{f(\bar{x} + \tau \bar d) - f(\bar{x})}{\tau}.$$
Given a possibly nonsmooth function $f : \R^n \mapsto \R\cup\{+\infty\}$ and a point $\bar{x}$ where $f(\bar{x})$ is finite, $df(\bar{x};\bar{d})$ is defined via
    $$df(\bar{x};\bar{d}) = \liminf_{\tau \searrow 0,~d \rightarrow \bar{d}} \frac{f(\bar{x} + \tau d) - f(\bar{x})}{\tau}.$$
This directional derivative is also known as the Hadamard lower derivative \cite{DeRu1995}, or the semiderivative \cite{rockafellar-wets-1998}.  Directional derivatives are linked to the subdifferential through the following classical formula
    \begin{equation}\label{eq:dfdef}
    \partial f(\bar{x}) = \{ v ~:~ v^\top d \leq df(\bar{x};d) ~\mbox{for all}~d\} \quad \mbox{\cite[Ex 8.4]{rockafellar-wets-1998}}.
    \end{equation}
Thus, given all possible directional derivatives, one can recreate the subdifferential as the infinite intersection of halfspaces,
    \[\partial f(\bar{x}) = \bigcap_{d : \|d \| = 1} \{ v \in \R^n : v^\top d \leq df(\bar{x}; d) \}.\]
Using this to develop a numerical approach to construct the subdifferential is, in general, impractical.  However, if the subdifferential is a polytope, then it may be possible to reconstruct the exact subdifferential using a finite number of directional derivative evaluations.

In general we will consider two basic cases:
\begin{enumerate}
    \item[I-] $\partial f(\bar{x})$ is a polytope, and an upper bound $\nf$ on its number of vertices $\nv$ is known;
    \item[II-] $\partial f(\bar{x})$ is a polytope, but no upper bound on the number of vertices is available.
\end{enumerate}
Case I corresponds to the situation where $f$ is a finite max function and the number of sub-functions $\nf$ used in constructing $f$ is known.  Case II corresponds to the case where $f$ is a finite max function, but no information about the function is available.
We shall see that the algorithms for both cases are the same, but Case I provides the potential for early termination.

\subsection{Links to Geometric Probing}\label{Geoprobing}

Under our assumption that $\partial f(\bar{x})$ is a nonempty polytope,
 and in light of Equation \eqref{eq:dfdef}, we  reformulate the problem in the following abstract manner:
\begin{quote} {\em Working in $\R^n$, the goal is to find all vertices $X_v = \{ v^1, v^2, \ldots, v^{\nv}\}$ of a nonempty compact polytope $X$ using an oracle $\oracle : \R^n \rightarrow \R$ that returns the value $\oracle(d) = \displaystyle \max_{v \in X} v^\top d$.}
\end{quote}
In general, polytopes will be denoted by a capital letter variable ($X, P, S$) and their corresponding sets of vertices will be denoted using a subscript $v$ (e.g., $X_v, P_v, S_v$).  In $\R^2$, we shall denote the edges of polytope $X$ by using a subscript $e$: $X_e = X \setminus (\mathrm{int}(X) \cup X_v)$.
Linking to constructing subdifferentials is achieved by setting
	$X = \partial f(\bar{x})$ and
	$\oracle(d) = df(\bar{x};d) $ $= \displaystyle \max_{v \in X} v^\top d$.
The assumptions for our two cases can be written as
\begin{enumerate}
    \item[I-] $X$ is a nonempty polytope with $\nv$ vertices, and $\nf \geq \nv$ is given,
    \item[II-] $X$ is a nonempty polytope with $\nv$ vertices, and $\nf = \infty$ is given.
\end{enumerate}

The interest in this reformulation is that it almost exactly corresponds to what is known as the {\em Hyperplane Probing} problem in the field of Geometric Probing.
In Hyperplane Probing the goal is to determine the shape of a polyhedral set in $\R^2$ given a hyperplane probe which measures the first time and place that a hyperplane moving parallel to itself intersects the polyhedral set \cite{Dobkin86, Skiena89}.  In Hyperplane Probing, it is generally assumed that $0$ is in the interior of the set, although this is only for convenience and has no effect on algorithm design \cite{Dobkin86}.

Hyperplane Probing dates back to at least 1986 \cite{Dobkin86}, where it was shown to be the dual problem to {\em Finger Probing} (where the probe measures the point where a ray exits the polyhedral set.  Using this knowledge, it was proven that to fully determine a polyhedral set with $\nv$ vertices, $3\nv - 1$ probes are necessary, and $3\nv$ probes are sufficient \cite{Dobkin86}.  Some variants of Hyperplane Probing exist.  In 1986, Bernstein considered the case when the polyhedral set $X$ is one of a finite list of potential sets $X \in \{X^1, X^2, ... X^p\}$ \cite{Bernstein1986}.  In this case the number of probes can be reduced to $2\nv + p$, where $\nv$ is the number of vertices of $X$ and $p$ is the size of the potential list of polyhedral sets.  In another variant, a double hyperplane probe is considered, which provides both the first and last place that a hyperplane moving parallel to itself intersects the polyhedral set \cite{Li1988}.  This extra information allows the resulting algorithm to terminate after  $3\nv-2$ probes.

Our problem differs from Hyperplane Probing in two small ways.  
First, instead of providing the first time {\em and place} that a hyperplane moving parallel to itself intersects the polyhedral set, our assumptions provide an oracle that yields the first time and but does not give the place.  Interestingly, this reduction of information has very little impact on the algorithm or convergence.  Indeed, we find it is sufficient to use $3\nv + 1$ oracle calls (as opposed to $3\nv$ for Hyperplane Probing).  In our case, the extra call is required to confirm that the final suspected vertex is indeed a vertex.  Second, we consider the space of polytopes $\R^n$, instead of polyhedral sets in $\R^2$.  While some recent research has examined Hyperplane Probing in $\R^3$ \cite{Imiya2012}, to our knowledge no research has explored the most general case of $\R^n$.
It is worth noting that the original work of Dobkin, Edelsbrunner, and Yap~\cite{Dobkin86} defined Hyperplane Probing as a problem in $\R^n$, but only studied the problem in $\R^2$.

\subsection{Notation}\label{notation}

Using the oracle notation of Geometric Probing, we introduce the following notation. The vector $\ej \in \R^n$ denotes the unit vector in the direction of the $j^{th}$ coordinate.  For a vector $d \in \R^n$ and a value $\oracle(d)$, we define the {\em generated constraint halfspace} by
    \[H(d) = \{ v \in \R^n ~:~ v^\top d \leq \oracle(d)\}\]
and the {\em generated hyperplane} by
    \[L(d) = \{ v \in \R^n ~:~ v^\top d = \oracle(d)\}.\]
Finally, given a set of vectors $D = \{d^1, d^2, ... d^m\}$ and corresponding values, we define the {\em generated constraint set} by
    \[\H(D) =\bigcap_{d \in D} H(d).\]
%
%

\section{One and two-dimensional  spaces}
\subsection{One-dimensional  space}\label{OneDim}

When working in $\R^1$ the problem is trivially solved.
Indeed, in $\R^1$, the number of vertices of $X$ must be either $1$ or $2$, i.e., the polytope $X$ will either be a single point, or a closed interval.

If $\nv$ is known to be equal to $1$, then a single evaluation suffices: $X = \{\oracle(1)\}$. If $\nv$ is unknown or is known to be $2$, then exactly two evaluations suffice. Specifically, evaluate $\oracle(1)$ and $\oracle(-1)$; if both are equal then $\nv=1$ and $X= \{\oracle(1)\}$, otherwise $\nv=2$ and $X = \conv\{\oracle(1), \oracle(-1)\}$.

\subsection{Two-dimensional  space}\label{TwoDim}

In $\R^2$, the problem becomes more difficult, as $\nv$ or its upper bound $\nf$ can take on any positive integer value.  
Our proposed algorithm (Algorithm \ref{Algo-exact} below) continually refines two polyhedral approximations of $X$.  The set $P$ is a polyhedral outer approximation of $X$, and the set $S$ is a polyhedral inner approximation of $X$.  The outer approximation set is initialized with a triangle, and the inner approximation is initialized as the empty set.  The algorithm proceeds by carefully truncating vertices of $P$ until a vertex of $P$ is proven to be a true vertex of $X$. As vertices of $P$ are found to be true vertices of $X$, they are added to $S$.
The algorithm terminates when $P = S$ or when the cardinality of $S_v$ equals $\nf$.

In the algorithm below, recall we denote the edge set of $P$ by $P_e$ and the vertex set of $P$ by $P_v$.  Similarly, $S_e$ is the edge set of $S$ and $S_v$ is the vertex set of $S$.

\begin{algo}{
	\label{Algo-exact}
	\sf Finding $X_v$ in $\R^2$, \\
    given an oracle $\oracle(d) = \displaystyle \max_{v \in X} v^\top d$, and \\
    given an upper bound on the number of vertices: $\nf \in \mathbb{N} \cup \{\infty\}$, $\nf \geq \nv$.}
Initialize: \\
\hspace*{5mm} \begin{tabular}[t]{|l}
 	Define the initial outer approximation polytope \\\hspace*{1cm}
	$P = \H(D)$ with $D  = \{ \e_1, \ \e_2, \ -\e_1-\e_2\}$.\\
    If $P$ is a singleton, then set $S_v = P$ and terminate.\\
    Otherwise, determine the 3 vertices of $P$ and enumerate them clockwise \\\hspace*{1cm}
    $P_v = \{p^1, p^2, p^3\}$.\\
    Create the (empty) initial inner approximation, and initialize counter \\
    \hspace*{1cm}$
    S_v = \emptyset$ and set $i=1$.\\
 \end{tabular}
	
{\sf While $P_v \neq S_v$ and $|S_v| < \nf$}\\
\hspace*{5mm} \begin{tabular}[t]{|l}
    Set $a = p^i$, $b=p^{i+1}$, and $c=p^{i+2}$, working modularly if $i+2 > |P_v|$.\\
	
	Choose $d \in \R^2$ such that $d^\top a = d^\top c  < d^\top b$. \\

	Compute $\oracle(d) \in [d^\top a, d^\top b]$.\\

	$\bullet$ If $\oracle(d) = d^\top b$, then $b \in X_v$:\\
	\hspace*{1cm} update $S_v \leftarrow S_v \cup \{ b \}$, $P_v \leftarrow P_v$,\\
    \hspace*{1cm} $i \leftarrow i+1$.\\

	$\bullet$ If $\oracle(d) = d^\top a$, then $\{a, c\} \in X_v$: \\
	\hspace*{1cm} update $S_v \leftarrow S_v \cup \{a, c\}$, $P_v \leftarrow P_v \setminus\{b\}$,\\
    \hspace*{1cm} re-enumerate $P_v$ starting at $a$ working clockwise,\\
    \hspace*{1cm} $i \leftarrow i+1$.\\
						
	$\bullet$ Otherwise (if $d^\top a = d^\top c < \oracle(d) <d^\top b$), then two new potential \\
	\hspace*{3mm}	vertices are discovered:\\
    \hspace*{1cm} compute the two distinct intersection points $\{b', c'\} = P_e \cap L(d)$,\\
	\hspace*{1cm} update $S_v \leftarrow S_v$, $P_v \leftarrow (P_v \setminus\{b\})\cup\{b', c'\}$,\\
    \hspace*{1cm} re-enumerate $P_v$ starting at $a$ working clockwise,\\
    \hspace*{1cm} $i \leftarrow i$.
\end{tabular}\\
Return $S_v$.
\end{algo}

Before examining the algorithm's convergence properties, we provide an illustrative example.

\clearpage
\newpage

\begin{Ex}\label{example}

Consider the polyhedral set $X$ in Figure \ref{fig:baseset}, a triangle in $\R^2$.

\begin{figure}[ht]
\begin{center}
\includegraphics[trim = 0mm 1mm 0mm 1mm, width=8cm]{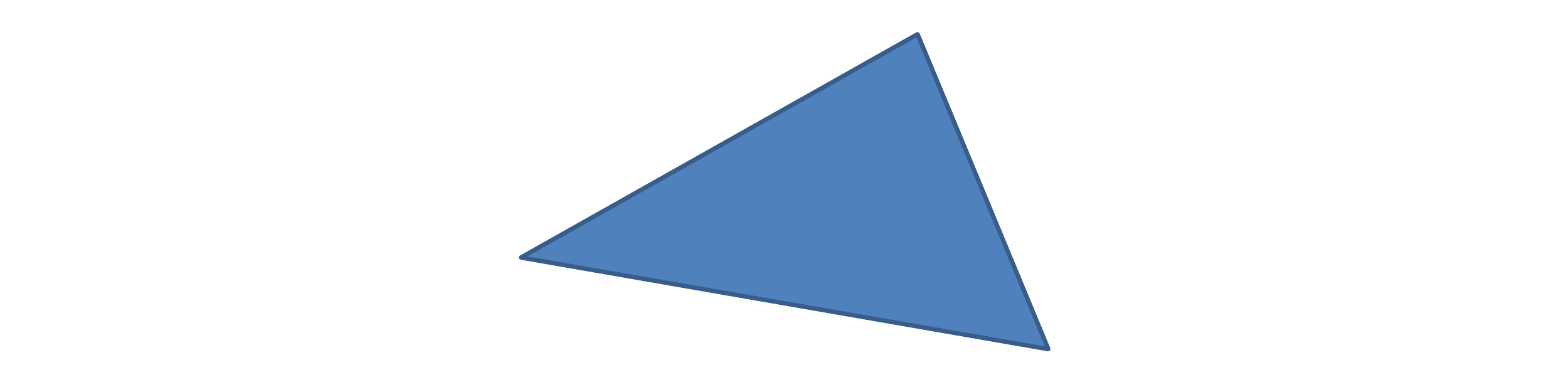}\end{center}
\caption{A polyhedral set $X$ in $\R^2$ used for Example \ref{example}.}\label{fig:baseset}
\end{figure}	

Using the notation of Subsection \ref{notation}, $P$ is initialized as $P = \H(D)$ with $D  = \{ \e_1, \ \e_2, \ -\e_1-\e_2\}$.  This creates a triangle containing $X$. If the triangle was degenerate, i.e., $P$ was a singleton, then the problem would be solved: $X = P$. In this example, $P$ is not degenerate, hence $P$ has exactly $3$ extreme points, which we label $\{p^1, p^2, p^3\}$, see Figure \ref{fig:initialize} (left). Notice that $P$ would still be a triangle if the set $X$ were a line segment rather than a triangle.

\begin{figure}[ht]
\begin{center}
	\includegraphics[trim = 30mm 1mm 5mm 1mm, clip, width=.48\textwidth]{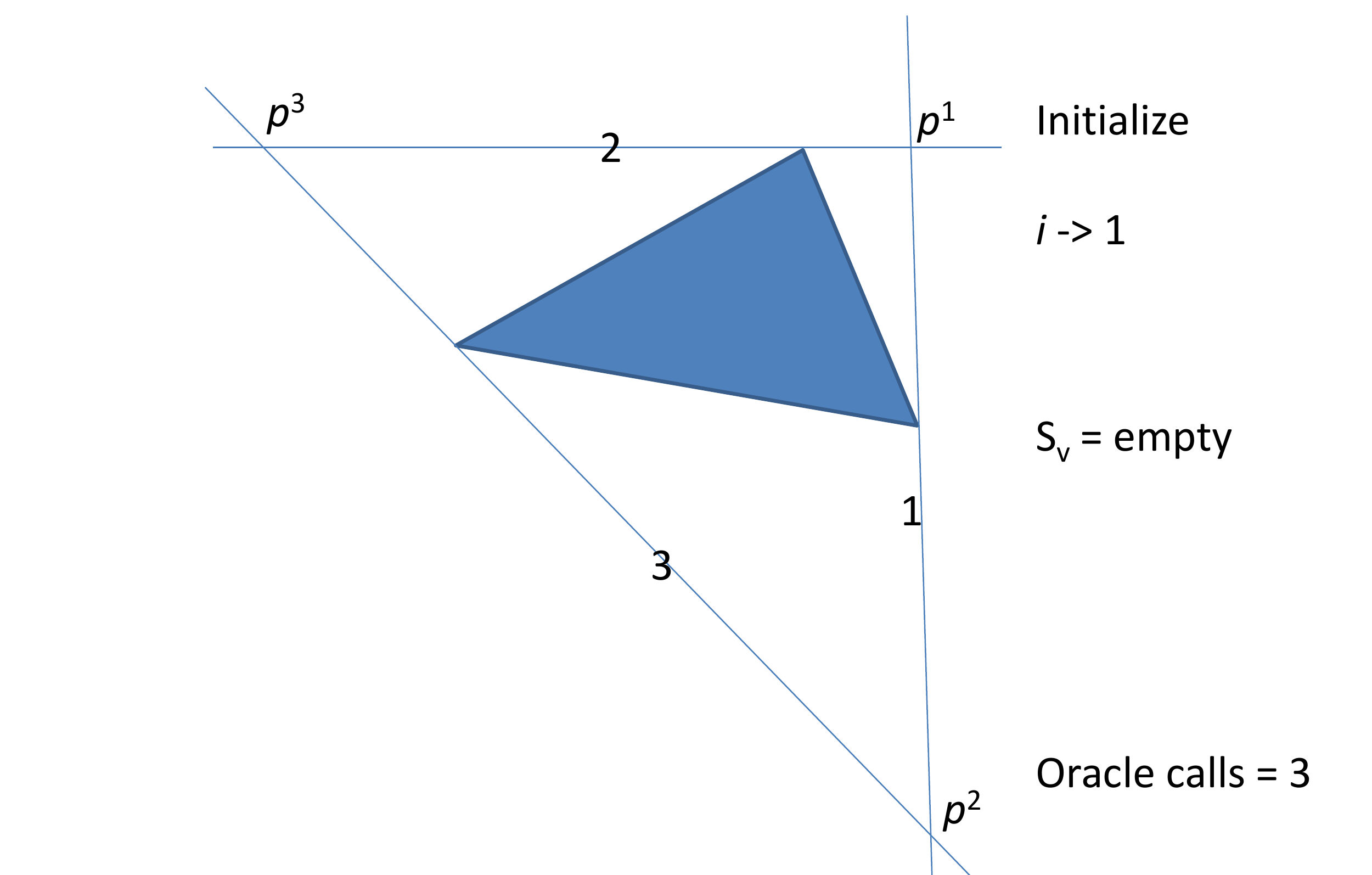} ~
	\includegraphics[trim = 30mm 1mm 5mm 1mm, clip, width=.48\textwidth]{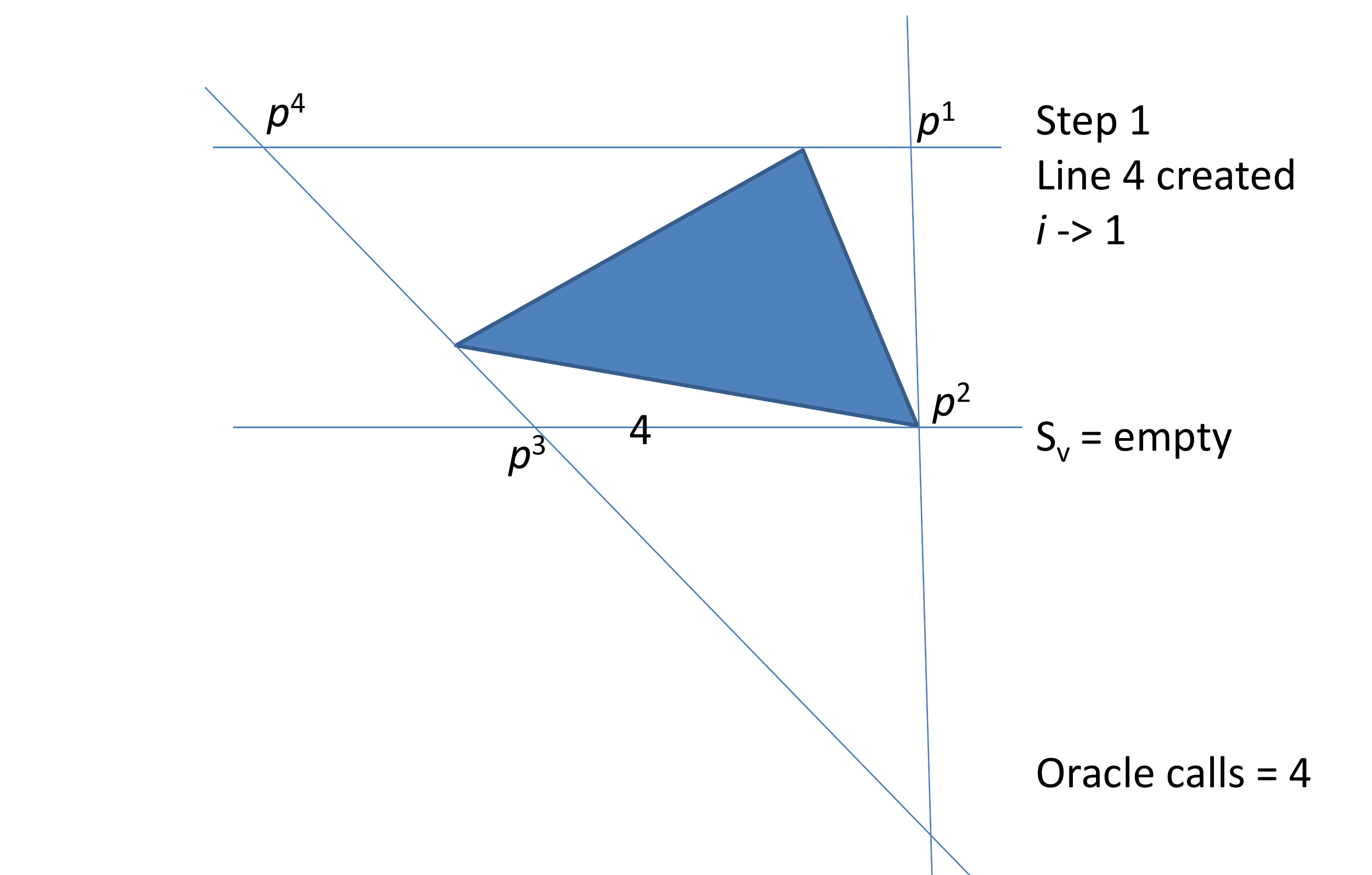}
\end{center}
\caption{[Left] Initialization of Algorithm \ref{Algo-exact}.
Numbers on the lines represent which oracle call created the tangent plane.
[Right] Step 1 of Algorithm \ref{Algo-exact}.
The new tangent plane is labelled with a 4 (the fourth oracle call).
}
\label{fig:initialize}
\end{figure}

Setting $i=1$, the algorithm creates a hyperplane parallel to the line segment adjoining $p^1$ and $p^3$ of Figure \ref{fig:initialize}: $d$ is set to $-\e_2$. The fourth oracle call is used to check for new potential vertices.  In this case $d^\top a < \oracle(d) < d^\top b$, so two new potential vertices are discovered. As a result, the index  $i$ is unchanged and only $P_v$ is updated.  The result appears in Figure \ref{fig:initialize} (right)
 in which the indices of the vertices are relabelled.


The reordered vertices in $P_v$ are now used to create a new hyperplane parallel to the line segment adjoining $p^1$ and $p^3$ (note $p^3$ has been relabelled from the left and right parts of Figure~\ref{fig:initialize}, so this line is distinct from the line in step 1).  In this case the result is $\oracle(d) = d^\top b$, where $b=p^2$.  As a result, the index $i$ is incremented to $2$, and $p^2$ is the first vertex to be discovered.  It is placed into $S_v = \{p^2\}$.  The result appears in Figure \ref{fig:step2} (left).

\clearpage \newpage

\begin{figure}[ht]
\begin{center}
	\includegraphics[trim = 30mm 1mm 5mm 1mm, clip, width=.48\textwidth]{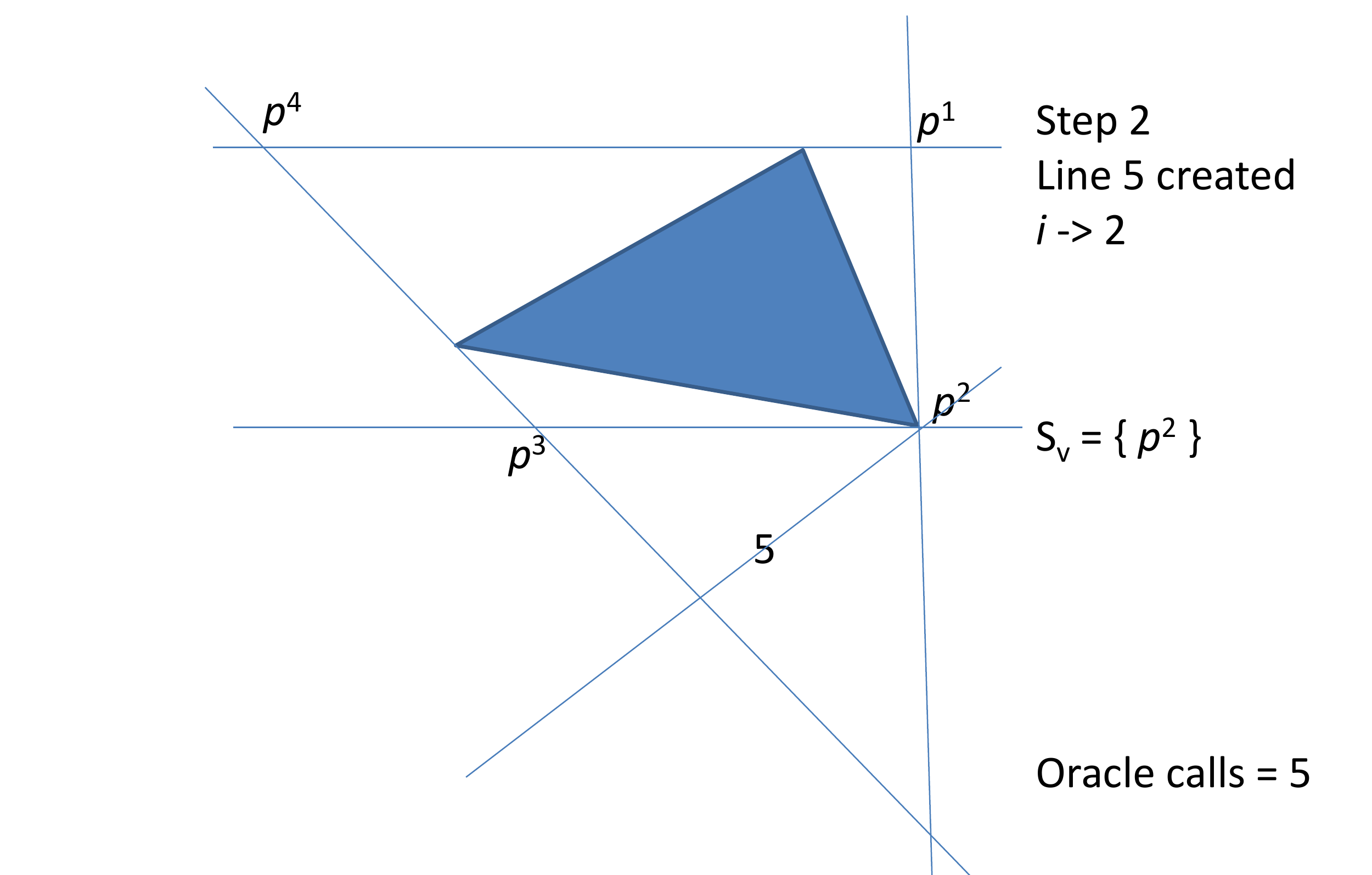} ~
	\includegraphics[trim = 30mm 1mm 5mm 1mm, clip, width=.48\textwidth]{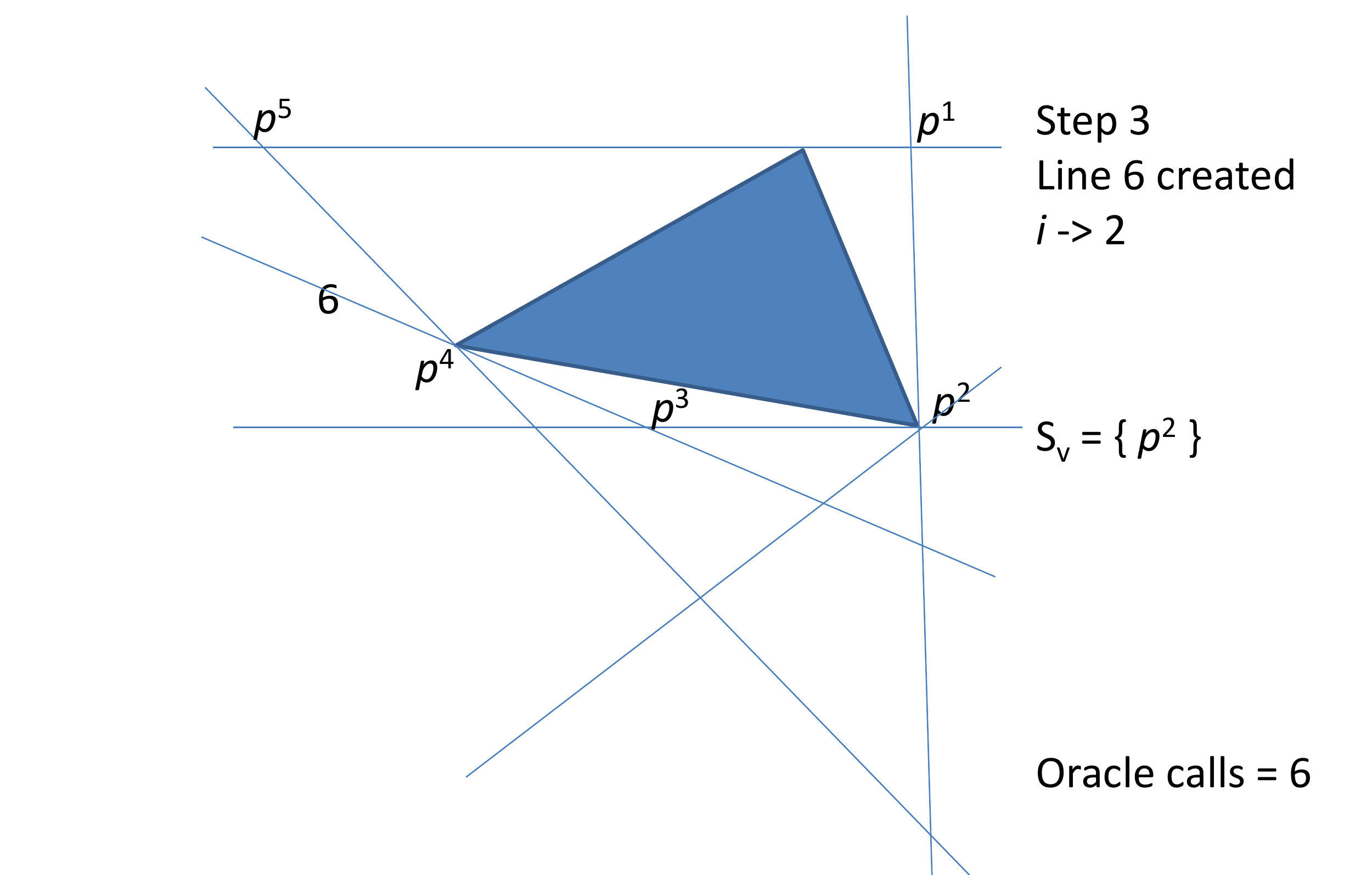}
\end{center}
\caption{[Left] Step 2 of Algorithm \ref{Algo-exact}.
The new tangent plane is labelled with a 5 (the fifth oracle call).
[Right] Step 3 of Algorithm \ref{Algo-exact}.
The new tangent plane is labelled with a 6 (the sixth oracle call).
}\label{fig:step2}
\end{figure}

In step 3, the new hyperplane is parallel to the segment joining $p^2$ and $p^4$ of Figure \ref{fig:step2} (left).
Like step 1, the result is two new potential vertices, and $i$ remains unchanged,
	see Figure \ref{fig:step2} (right).


Step 4, Figure \ref{fig:step4} (left),
	discovers vertex $p^3$ with an hyperplane tangent to one of the sides of the triangle.
The index $i$ is incremented, and $p^3$ is added to $S_v$.

\begin{figure}[ht]
\begin{center}
	\includegraphics[trim = 30mm 1mm 5mm 1mm, clip, width=.48\textwidth]{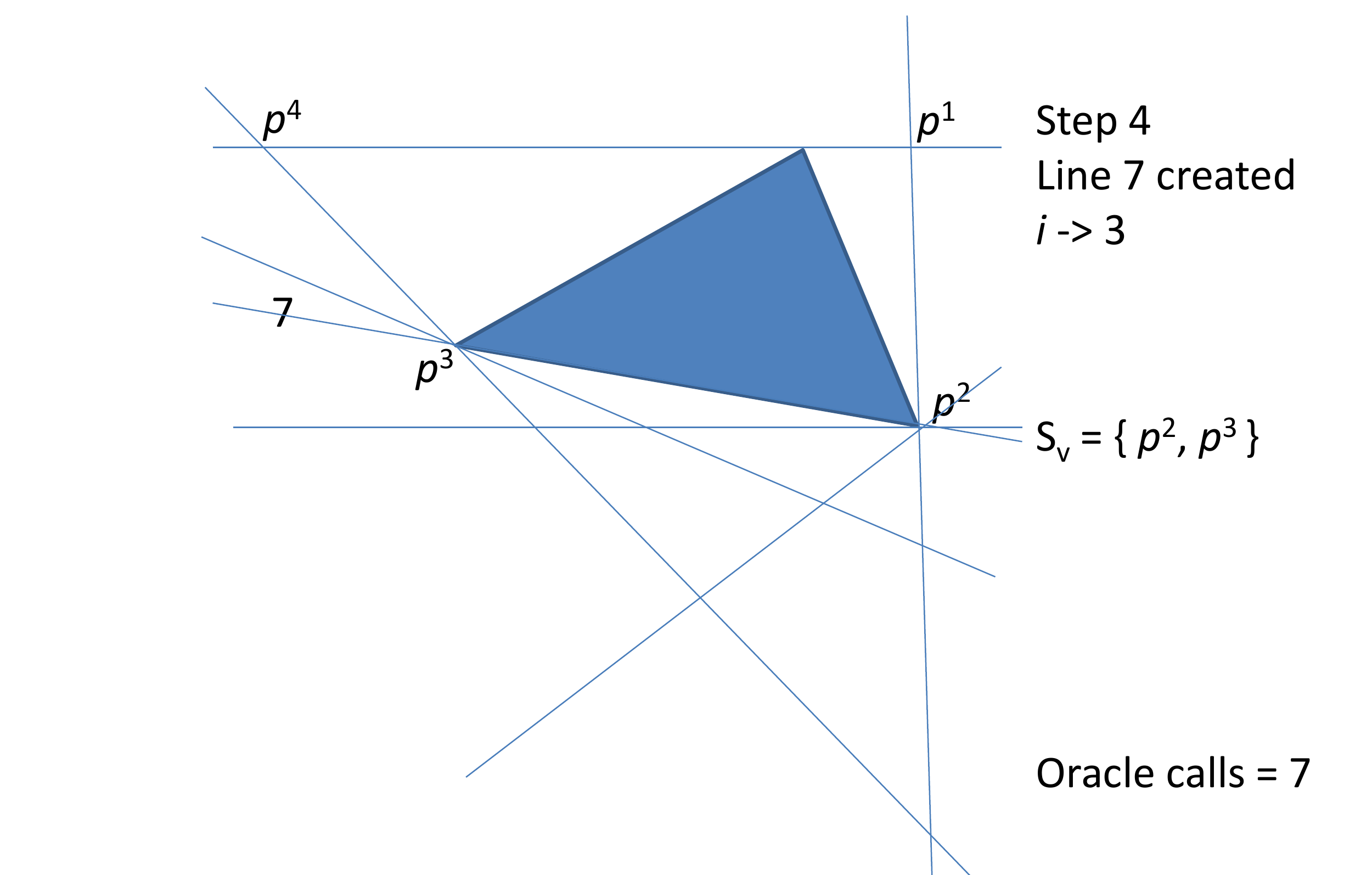} ~
	\includegraphics[trim = 30mm 1mm 5mm 1mm, clip, width=.48\textwidth]{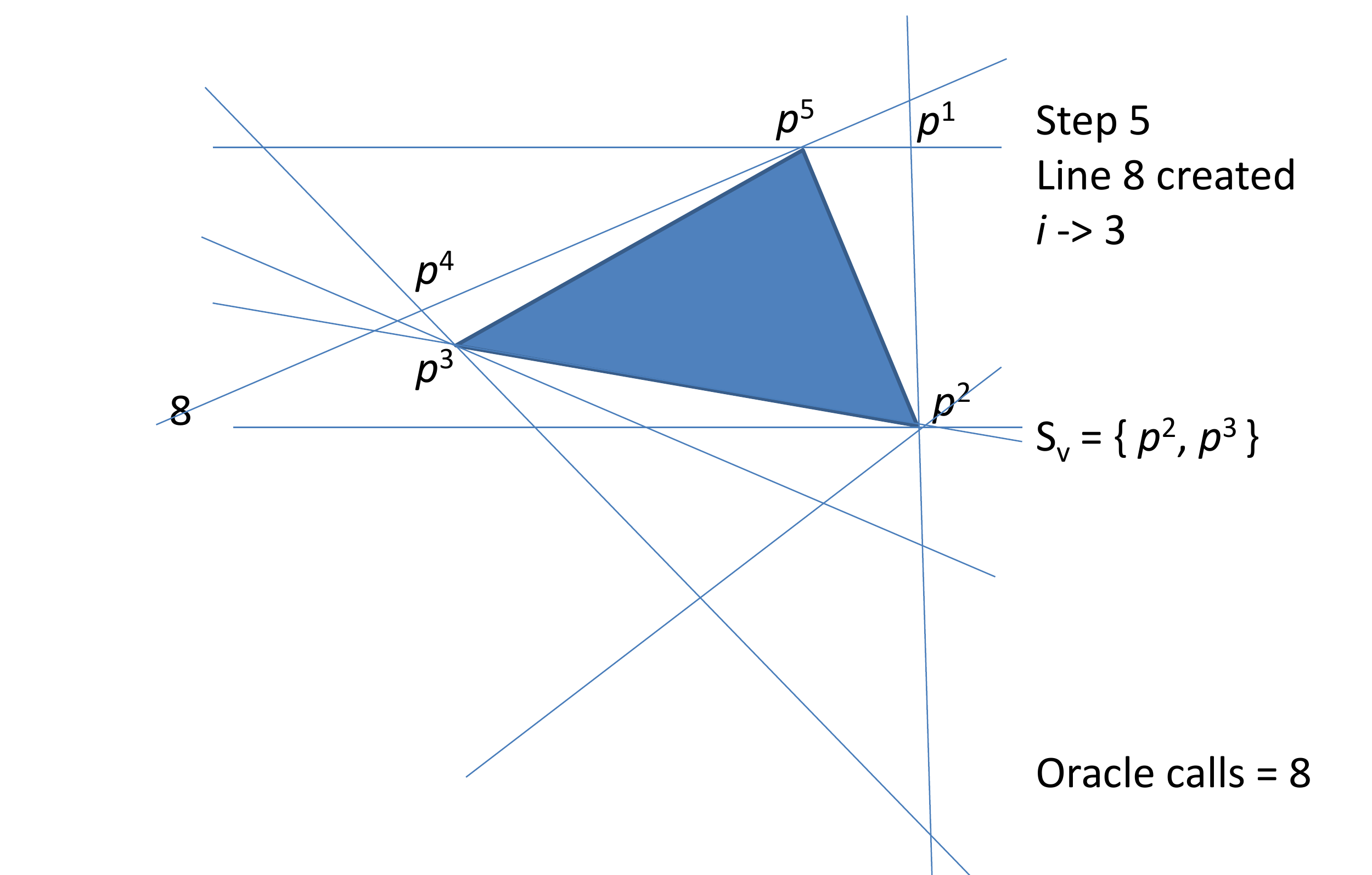}
\end{center}
\caption{[Left] Step 4 of Algorithm \ref{Algo-exact}.
The new tangent plane is labelled with a 7 (the seventh oracle call).
[Right] Step 5 of Algorithm \ref{Algo-exact}.
The new tangent plane is labelled with a 8 (the eigth oracle call).
}\label{fig:step4}
\end{figure}


Step 5, Figure \ref{fig:step4} (right), creates two new potential vertices with an hyperplane parallel to the line segment joining the $p^3$ an $p^1$ from step 4.
This is an example of using modularity.


In step 6, Figure \ref{fig:step6} (left), we have $d^\top a = d^\top c = \oracle(d)$, so two vertices are discovered simultaneously.
(In this case, vertex $p^3$ was actually already known, but this is not necessarily always the case.)
This results in $p^5$ being relabelled as $p^4$ and being added to $S_v$.

\begin{figure}[ht]
\begin{center}
	\includegraphics[trim = 30mm 1mm 5mm 1mm, clip, width=.48\textwidth]{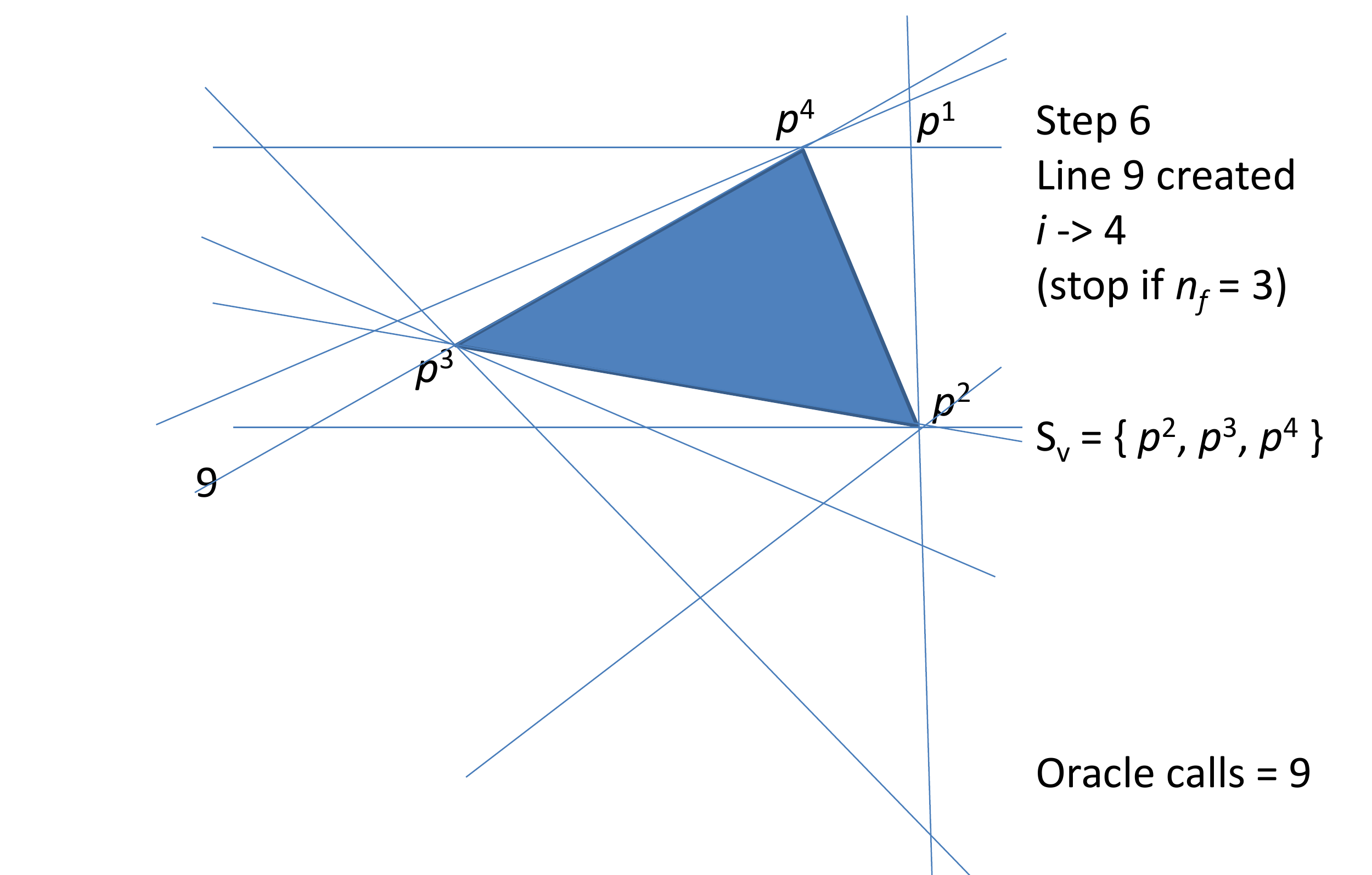} ~
	\includegraphics[trim = 30mm 1mm 5mm 1mm, clip, width=.48\textwidth]{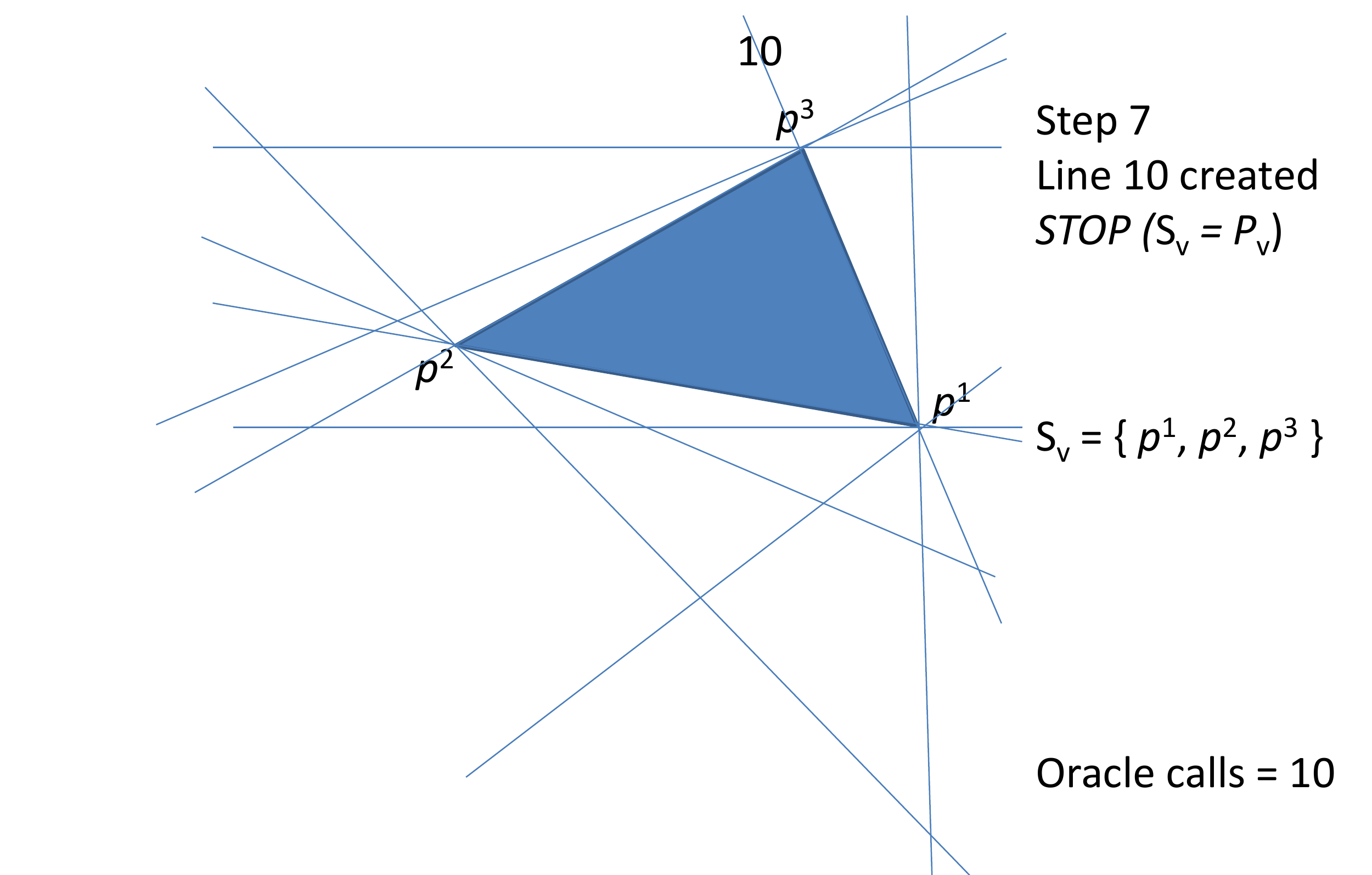}
\end{center}
\caption{[Left] Step 6 of Algorithm \ref{Algo-exact}.
The new tangent plane is labelled with a 9 (the ninth oracle call).
[Right] Step 7 of Algorithm \ref{Algo-exact}.
The new tangent plane is labelled with a 10 (the tenth oracle call).
}\label{fig:step6}
\end{figure}

If $\nf = 3$, then the algorithm can stop at this point, as $S_v$ contains three known vertices.  Conversely, if $\nf > 3$, then the algorithm requires one more step to truncate  $p^1$ from the potential vertex list, see Figure \ref{fig:step6} (right).


In summary, the algorithm uses either a total of $3\nv=9$ oracle calls if $\nf=3$,
	or $3\nv+1 = 10$ oracle calls if $\nf > 3$.
In both situations it returns $S_v = X_v$.
\end{Ex}

Example \ref{example} demonstrates the ideas behind the algorithm, and shows that it is possible to require $3\nv+1$ oracle calls. The next theorem proves the algorithm converges to the correct vertex set.  It also proves that, if $\nf=\nv$, then at most $3\nv$ oracle calls will be required, and if $\nf > \nv$, then at most $3\nv+1$ oracle calls will be required.

\begin{Th}[Convergence of Algorithm \ref{Algo-exact}] \label{Th:exactbounds}
Let $X$ be a polytope in $\R^2$ with $\nv$ vertices contained in set $X_v$.

\noindent
If $\nv=1$, then Algorithm~\ref{Algo-exact} terminates after the initialization step.
\\
If $\nf=\nv$, then Algorithm~\ref{Algo-exact} terminates after at most $3\nv$ oracle calls.
\\
If $\nf > \nv$, then Algorithm~\ref{Algo-exact} terminates after at most $3\nv+1$ oracle calls.
\\
In either case,  the algorithm terminates with $S_v = X_v$.
\end{Th}

\proof We shall use the notation of Subsection \ref{notation}.
First note, if $X_v$ is a singleton, then the initialization step will result in $P = X_v$ and the algorithm terminates  after $3$ oracle calls.

If $X_v$ is not a singleton, then each oracle call of the algorithm, $\oracle(d)$ introduces a new tangent plane to $X$.  Specifically
    $$L(d) = \{ v \in \R^2 ~:~ v^\top d = \oracle(d) \}$$
is a tangent plane to $X$. As $X$ is polyhedral, we must have $L(d) \cap X_v \neq \emptyset$. Let $v \in L(d) \cap X_v$. The vertex $v \in X_v$ lies in one of three sets: the interior, the edges or the vertices of $P$.

If $v \in \mathrm{int} P$, then the previously undiscovered vertex $v$ of $X_v$ has been added to $P_e$. As $X$ has $\nv$ vertices, this can happen at most $\nv$ times.

If $v \in P_e$, then it will be shifted from $P_e$ into the potential vertex set $P_v$. Again, as $X$ has $\nv$ vertices, this can happen at most $\nv$ times.

Finally, if $v \in P_v$, then we are in the case of $D(d) = d^\top b$ or $D(d) = d^\top a$, so the potential vertex $v \in P_v$ has been confirmed as a true vertex of $X$ and placed in $S_v$. This can happen at most $\nv$ times.

Thus, after at most $3 \nv$ oracle calls, $S_v$ will contain all $\nv$ vertices of $X$. If $\nf=\nv$, then the algorithm will terminate at this point.

If $\nf > \nv$ and $p^1 \notin X_v$, then after at most $3\nv$ oracle calls, $S_v=X_v$, but $P_v$ may still contain $p^1$. One final oracle call will remove $p^1$ from $P_v$,  making $S_v = P_v$ and the algorithm will terminate.\qed

In some situations it is possible to terminate the algorithm early.

\begin{Lem}[Improved stopping when $\nf = \nv$] \label{Lem:nfequalsnv}
Let $X$ be a polytope with $\nv$ vertices contained in set $X_v$.  Suppose $\nf = \nv$.  Suppose the algorithm has run to the point where $\nv-1$ vertices are identified.  If $P_e$ contains two edges that are not adjacent to any of the known vertices, then the intersection of those two edges must be the final vertex.  As such, the algorithm can be terminated.
\end{Lem}

Lemma \ref{Lem:nfequalsnv} is particularly useful when $\nf = \nv = 2$.

\begin{Cor}[Special case of $\nf = \nv = 2$]\label{Cor:nf2}Let $X$ be a polytope with $\nv = 2$ vertices contained in set $X_v$. If $\nf=2$, then the algorithm can be terminated after just $5$ oracle calls.
\end{Cor}

\begin{figure}[ht]
\begin{center}
\includegraphics[height=2in]{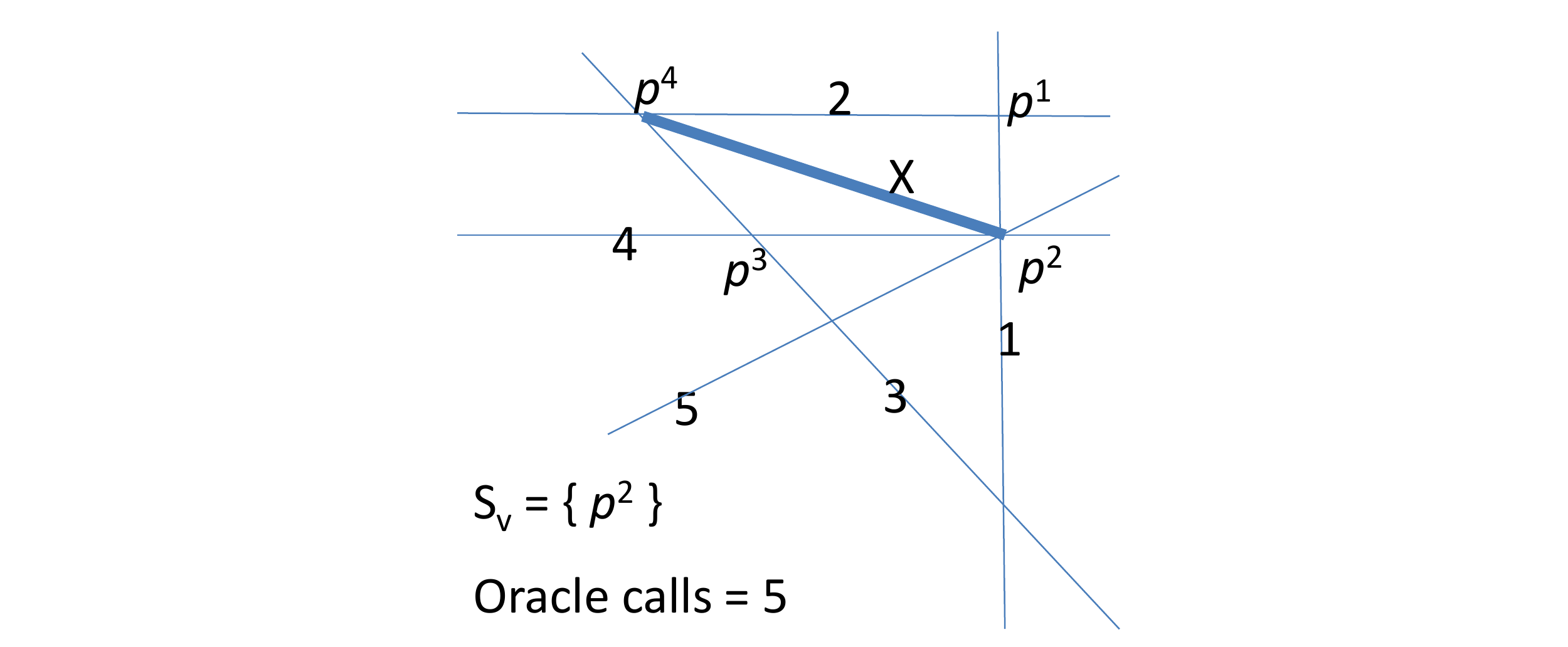}
\end{center}
\caption{The situation after $5$ oracle calls if $\nv=\nf=2$ and only one vertex has been identified.  By Lemma \ref{Lem:nfequalsnv}, potential vertex $p^4$ must be a true vertex.}\label{fig:5calls}
\end{figure}

\proof Following the logic in the proof of Theorem \ref{Th:exactbounds}, at most $2$ oracle calls can move a vertex from the interior of $P$ to the edge set of $P$, and at most $2$ oracle calls can move a vertex from the edge set of $P$ to the vertex set of $P$. Therefore, after 5 oracle calls, at least one vertex has been identified.  If two vertices are identified, then we are done.  Otherwise, we must be in the situation shown in Figure \ref{fig:5calls}.

\clearpage
\newpage

In particular, we must have one vertex of $P$ that has $3$ lines through it, one of which is redundant in defining $P$.  The other 2 lines that make up $P_e$ must not intersect this vertex, and cannot create a vertex of $P$ with 3 lines through it.  The only possible way to do this is a quadrilateral.  Lemma \ref{Lem:nfequalsnv} now applies, so the final vertex can be identified without an additional oracle calls. \qed

It is worth noting that, unless the special termination trick in Lemma \ref{Lem:nfequalsnv} applies, then the bounds provided in Theorem \ref{Th:exactbounds} are tight, as was demonstrated in Example \ref{example}.


It is clear that there is nothing particularly special about the directions $e_1, e_2,$ and $-e_1-e_2$ used in the initialization of Algorithm \ref{Algo-exact}.  If these directions are replaced by any three directions that {\em positively span} $\R^2$, then the algorithm behaviour is essentially unchanged. More interestingly, the initialization set can be replaced by any set that positively spans $\R^2$, with the only negative impact being the potential to waste oracle calls during the initialization phase.
The following lemma analyses this situation and will be referred to later in the paper.

\begin{Lem}\label{lem:changeinitializationM}
Suppose the initialization set $\{e_1, e_2, -e_1-e_2\}$ in Algorithm \ref{Algo-exact} is replaced by $D = \{d^1, d^2, ..., d^m\}$, where $D$ {\em positively spans} $\R^2$ in the sense that $$\{ x ~:~ x = \sum_{i=1}^m \lambda_i d^i, \lambda^i \geq 0\} = \R^2.$$

\noindent If $\nf=\nv$, then Algorithm~\ref{Algo-exact} terminates after at most $(3\nv) + (m-3)$ oracle calls.

\noindent If $\nf > \nv$, then Algorithm~\ref{Algo-exact} terminates after at most $(3\nv+1) + (m-3)$ oracle calls.

\noindent In either case, when the algorithm terminates we have $S_v = X_v$.
\end{Lem}

\proof  Since $\{ x ~:~ x = \sum_{i=1}^m \lambda_i d^i, \lambda^i \geq 0\} = \R^2$, the initialization directions allow the algorithm to create a compact initialization polytope $P$ that contains $X$.  If a vertex of $P$ is defined by $4$ or more oracle, then any oracle call past the first $3$ is potentially wasted.  However, all other oracle calls follow the same rules as those in the proof of Theorem \ref{Th:exactbounds}.  The bounds follow from the fact that a maximum of $m-3$ oracle calls will be wasted. \qed

We conclude this section with a remark that will be important in analyzing Algorithm \ref{Algo-Rn-nf3}.

\begin{Rem}\label{Rem:changeinitializationAxis}
If the initialization set $D = \{e_1, e_2, -e_1, -e_2\}$ is used and $\nv \geq 2$, then no oracle calls will be wasted in creating the initialization set.  Thus, in this case, the exact bounds of Theorem \ref{Th:exactbounds} will still hold.  (If $\nv =1$, then the resulting will use $4$ oracle calls instead of $3$, but the algorithm will still terminate immediately after the initialization step.)
\end{Rem}

%

\section{$n$-dimensional  space with at most $3$ vertices}\label{nDim}
This section is devoted to $\R^n$ in which the polytope $X$ has at most $\nf \leq 3$ vertices.  We consider $3$ subcases.

\subsection{$n$-dimensional  subcase with $\nf=1$}\label{nDim-nf1}
The simplest case in $\R^n$ occurs  when the upper bound on the number of vertices of $X$ is $\nf = 1$.
This trivial case is solved by $n$ calls to the oracle $\oracle$:
For $i = 1, 2, ... n$, evaluate $\oracle(\e_i)$ and return $X_v = \{\oracle(\e_1)\e_1 + \oracle(\e_2)\e_2 + ... + \oracle(\e_n)\e_n\}.$

\subsection{$n$-dimensional  subcase with $\nf=2$}\label{nDim-nf2}

The next simplest case is when the upper bound is $\nf = 2$.  We present an alternate algorithm for this case, which uses a similar approach to constructs an outer approximation $P$ of $X$.  However, in this case a hyperrectangle is used in the initialization phase to bound the vertices.

\begin{algo}{
	\label{Algo-Rn-nf2}
    \sf Finding $X_v$ in $\R^n$, \\
     given an oracle $\oracle(d) = \displaystyle \max_{v \in X} v^\top d$, and \\
     given the upper bound on the number of vertices is $\nf = 2 \geq \nv$.}

Initialize: \\
\hspace*{5mm} \begin{tabular}[t]{|l}
 	Define the initial outer approximation polytope hyperrectangle \\\hspace*{1cm}
	$P = \H(D)$ with $D  = \{ \pm\e_i : i = 1,2\ldots, n\}$.\\
    If $P$ is a singleton, then return $P$ and terminate.\\
    Otherwise, set $\ell_i = -\oracle(-\e_i)$ and $u_i = \oracle(\e_i)$ for all $i$.\\
    By relabelling indices if necessary, assume $\ell_1 < u_1$.\\
    Initialize points $a = (\ell_1, \ell_2, ..., \ell_n)$, $b = (u_1, u_2, ..., u_n)$.
 \end{tabular}
	
For $i$ from $2$ to $n$\\
\hspace*{5mm} \begin{tabular}[t]{|l}
	If $\ell_i < u_i$, then\\
	\hspace*{1cm}
		choose $d \in \R^n$ with $d_i = 1$ and $d_j=0$ for all indices $j > i$ \\
	\hspace*{1cm}
		such that 
		$d^\top a = d^\top b$; \\
	
	\hspace*{1cm}
		if $\oracle(d) > d^\top a$, then reset
		$a_i \leftarrow u_i$ and $b_i \leftarrow \ell_i$.\\
\end{tabular}\\
Return $\{a, b\}$.
\end{algo}

\begin{Th}[Convergence of Algorithm \ref{Algo-Rn-nf2}] \label{Th-Rn-nf2}
Let $X$ be a polytope with $\nv \in \{1, 2\}$ vertices.  Let the known upper bound for $\nv$ be $\nf=2$.

\noindent
If $\nv=1$, then Algorithm~\ref{Algo-Rn-nf2} terminates after exactly $2n$ oracle calls with the singleton $X_v = P$.
\\
If $\nv=2$, then Algorithm~\ref{Algo-Rn-nf2} terminates after at most $3n-1$ oracle calls with $X_v = \{ a, b \}$.

\end{Th}

\proof The initialization phase of the algorithm calls the oracle exactly $2n$ times, providing bounds $\ell_i \leq x_i \leq u_i$ for each index $i$.  It is obvious that $X$ is a singleton if and only if $\ell_i=u_i$ for all $i \in \{1, 2, ..., n\}$, in which case $X=P$. Otherwise the algorithm proceeds into the iterative loop, that identifies the components of the two vertices of $X_v$.

For each $i \in \{2, ... n\}$, the bound $\ell_i$ is initially assigned to $a_i$, and $u_i$ is initially assigned to $b_i$. If both lower and upper bounds are equal, $\ell_i = u_i$, then both $a_i$ and $b_i$ remain at this value.  When they differ, an additional call to the oracle will indicate to which of $a_i$ or $b_i$ will the bounds be associated.  The algorithm constructs the vector $d \in \R^n$ in such a way that $d^\top a = d^\top b$, which implies that $\sum_{k=1}^{i-1} d_i a_i + \ell_i = \sum_{k=1}^{i-1} d_i b_i + u_i$.  Therefore, if $\oracle(d) = d^\top a$, then it follows that the lower bound was correctly attributed to $a_i$ and the upper bound was correctly attributed to $b_i$. Otherwise, they need to be swapped: $a_i \leftarrow u_i$ and $b_i \leftarrow \ell_i$.

The total number of calls to the oracle is $2n$ for initialization step plus a maximum of $n-1$ from the loop of $i = 2$ to $n$.  Thus, an overall  maximum of $3n-1$ oracle calls are required. \qed


\subsection{$n$-dimensional  subcase with $\nf=3$}\label{nDim-nf3}

Reconstructing the subdifferential gets more complicated as the number of vertices increases.
In this last subcase, we propose a method for the $n$-dimensional case
	where the number of vertices is bounded by $\nf=3$.
The proposed strategy exploits the following fact.
Let $Y$ be the projection on a linear subspace of the polytope $X \subset \R^n$.
Any vertex $y$ of $Y$ is the projection of a vertex $x$ of $X$~\cite{Tiwary2008}.
If the number  of vertices is small ($\nv \in \{1, 2, 3\}$), then the combinatorics involved in deducing the vertices of $X$ from those of $Y$ is manageable.

We propose the following algorithm
 that proceeds by successively finding the vertices of the projections of $X$ on
 $\R^2, \R^3 \ldots, \R^n$.
The number of vertices of the projection in $\R^k$ is denoted $n^k_v$.

\begin{algo}{
	\label{Algo-Rn-nf3}
    \sf Finding $X_v$ in $\R^n$, \\
     given an oracle $\oracle(d) = \displaystyle \max_{v \in X} v^\top d$, and \\
     given the upper bound on the number of vertices is $\nf = 3 \geq \nv$.}

Initialize: \\
\hspace*{5mm} \begin{tabular}[t]{|l}
	Apply Algorithm~\ref{Algo-exact} to obtain the $n^2_v$
	vertices of the projection of $X$ in $\R^2$:\\
	$x^j = (x^j_1, x^j_2, 0,0,\ldots,0) \in \R^n$ for $1 \leq j \leq n^2_v$. \\
	\end{tabular}
	
For $k$ from $3$ to $n$\\
\hspace*{5mm} \begin{tabular}[t]{|l}
	Given $x^j = (x^j_1, x^j_2, \ldots, x^j_{k-1} 0,0,\ldots,0) \in \R^n$ for $1 \leq j \leq n^{k-1}_v$. \\
	Call the oracle twice and set
		$\ell_k=  \oracle(-\e_k)$ and $u_k= \oracle(\e_k)$.\\
	If $\ell_k = u_k$, then\\
	\hspace*{5mm} \begin{tabular}[t]{l}
		set $n^k_v = n^{k-1}_v$ and reset
			$x^j_k \leftarrow \ell_k$ for each $1 \leq j \leq n^k_v$.\\
	\end{tabular}\\
	Otherwise $\ell_k < u_k$ and apply the appropriate case\\
	\hspace*{5mm} \begin{tabular}[t]{l}
		\underline{$\bullet$  Case I : $n^{k-1}_v = 1$.} \\
		\hspace*{5mm} \begin{tabular}[t]{l}
			Set $n^k_v = 2$, $x^2 = x^1$
			and reset  $x^1_k \leftarrow \ell_k$ and $x^2_k \leftarrow u_k$.
			\end{tabular}\\
		\underline{$\bullet$ Case II : $n^{k-1}_v = 2$.}\\
		\hspace*{5mm} \begin{tabular}[t]{l}
			The $n^k_v \leq 3$ vertices of the projection in $\R^k$
			lie in the two-dimen-\\sional plane containing
			$\{ x^i + \ell_k \e_k, x^i + u_k \e_k \ : \ i =1,2 \}$.	\\	
			Use a change of variables to reduce this to a problem in $\R^2$ and \\
            apply Algorithm \ref{Algo-exact} using the current initialization state, as\\
              allowed by Lemma \ref{lem:changeinitializationM}.
			\end{tabular}\\
		\underline{$\bullet$ Case III : $n^{k-1}_v = 3$.} \\
        \hspace*{5mm} Set $n^k_v = 3$.\\
		\hspace*{5mm} For j = 1,2,3\\
		\hspace*{10mm} \begin{tabular}[t]{|l}	
			Choose  $d \in \R^n$ so that
			$d^{\top}(x^j + \ell_k \e_k) = 0,$ and \\
			$d^{\top}(x^i + u_k \e_k) = \left\{ \begin{array}{ll}
				1 & \mbox{ if } \ i =j \\
				0 & \mbox{ if } \ i \in \{1,2,3\}, \  i \ne j.
				\end{array} \right.$\\
			Call the oracle and reset $x^j_k \leftarrow \ell_k + \oracle(d)(u_k - \ell_k)$.
		\end{tabular}
	\end{tabular}\\
\end{tabular}\\
Return $S_v = \{x^j : 1 \leq j \leq n^n_v\}$.
\end{algo}

\begin{Th}[Convergence of Algorithm \ref{Algo-Rn-nf3}] \label{Th:Rn-nf3}
Let $X$ be a polytope in $\R^n$ with $\nv \leq \nf = 3$ vertices contained in set $X_v$.

\noindent
If $\nv = 1$, then Algorithm~\ref{Algo-Rn-nf3} terminates after at most $2n-1$ oracle calls.
\\
If $\nv = 2$, then Algorithm~\ref{Algo-Rn-nf3} terminates after at most $5n-3$ oracle calls.
\\
If $\nv = 3$, then Algorithm~\ref{Algo-Rn-nf3} terminates after at most $5n-1$ oracle calls.
\\
In all cases, the algorithm terminates with $X_v = S_v$.
\end{Th}

\proof The proof is done by induction on $n$, the dimension of the space.

If $n=2$, then the algorithm reduces to Algorithm \ref{Algo-exact} and the bounds follow from Theorem \ref{Th:exactbounds}.

Now, suppose that the result is true for some $n = k-1 \geq 2$.  Set
    $$
    \ZZ = \left\{\begin{array}{ll}
    2n-1 & \mbox{if}~\nv^{n}=1,\\
    5n-3 & \mbox{if}~\nv^{n}=2,\\
    5n-1 & \mbox{if}~\nv^{n}=3,\\
    \end{array}\right.
    = \left\{\begin{array}{ll}
    2k-3 & \mbox{if}~\nv^{k-1}=1,\\
    5k-8 & \mbox{if}~\nv^{k-1}=2,\\
    5k-6 & \mbox{if}~\nv^{k-1}=3.\\
    \end{array}\right.
    $$
Thus, after $\ZZ$ oracle calls, the algorithm has correctly identified the $\nv^{k-1}$ vertices of $X$ projected onto $\R^{n} = \R^{k-1}$ as a subspace of $\R^k$.

As the Algorithm~\ref{Algo-Rn-nf3} proceeds into the $k^{\mbox{\scriptsize th}}$ iterative loop, it uses $2$ oracle calls to compute the bounds  $\ell_k$ and $u_k$ on the $k^{\mbox{\scriptsize th}}$ variable.
If both values are equal, then this implies that all vertices belong to the subspace where $x_k = \ell_k$ and no more calls to the oracle are required, and the algorithm terminates using $\ZZ+2$ oracle calls.
If $\ell_k < u_k$, then the algorithm proceeds to Case I, II, or III.

In Case I, $\ZZ=2k-3$ oracle calls sufficed to find the $n^{k-1}_v$ vertices in the projection on $\R^{k-1}$.
No other calls are required because all vertices in $\R^k$ belong to the one dimensional subspace in which the $k^{\mbox{\scriptsize th}}$ components lie in the interval $[x_k + \ell_k, x_k + u_k]$.  The overall number of oracle calls is $\ZZ+2 = 2k-1$.

In Case II, $\ZZ=5k-8$, and the $n^k_v$ vertices of the projection in $\R^k$
 must lie in the two dimensional plane containing $\{ x^1 + \ell_k \e_k, x^1 + u_k \e_k, x^2 + \ell_k \e_k , x^2 + u_k \e_k\}$.
 Remark~\ref{Rem:changeinitializationAxis} ensures that at most $3$ additional  oracles calls are required if $\nv^k = 2$ and at most $5$ more oracles calls are required if $\nv^k = 3$.   So, the algorithm terminates using at most $\ZZ+5 = 5k-3$ oracle calls if $\nv=2$, and $\ZZ+7 = 5k-1$ oracle calls if $\nv=3$.

\begin{figure}[htb!]
\centering
	\includegraphics[width=8.7cm]{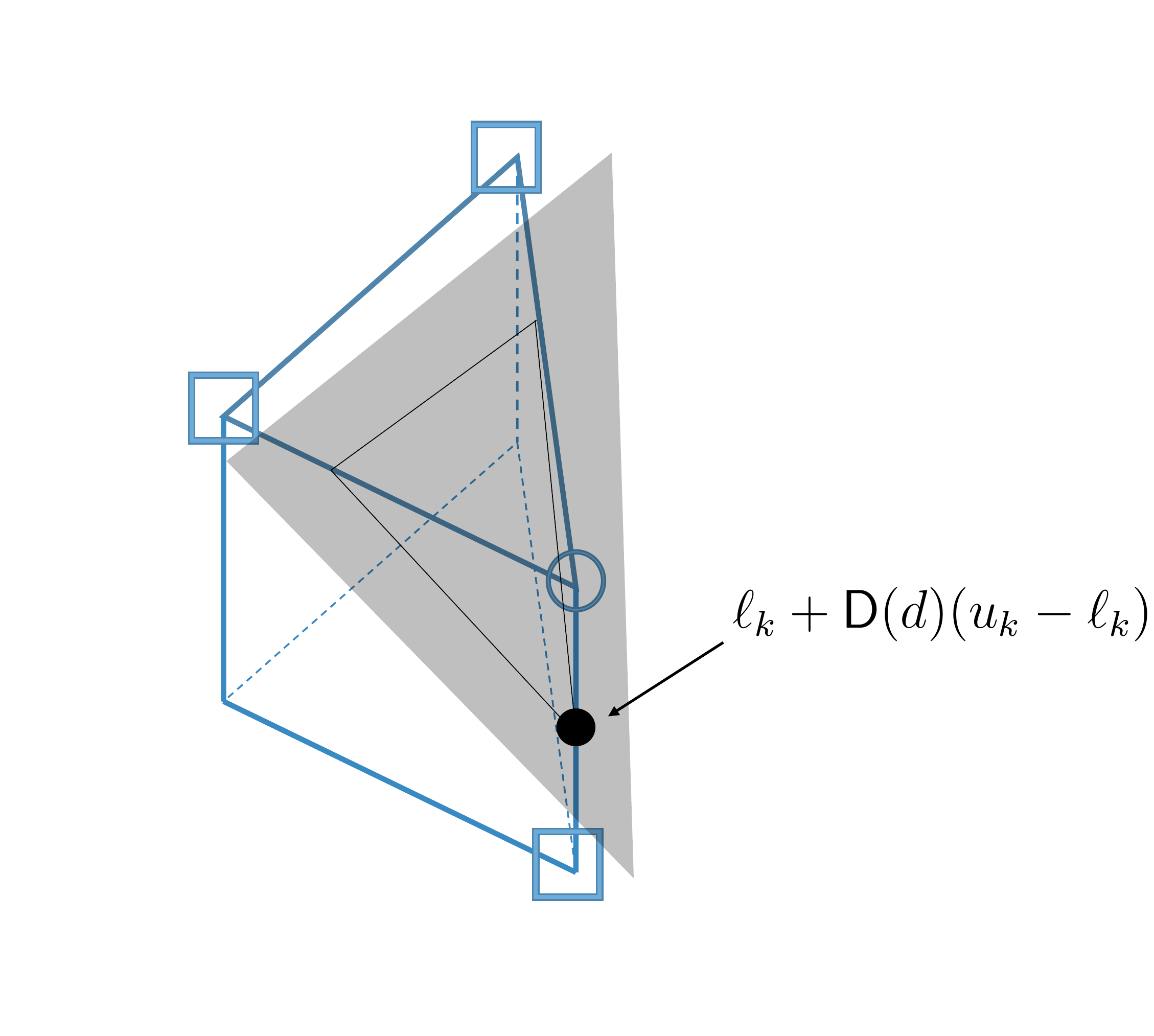}
\caption{Case III: all three vertices are located on the vertical edges of a triangular prism.  The triangular prism is the outer bound created by the initialization stage of algorithm \ref{Algo-Rn-nf3}.  The set $X$ is the triangle within the prism and the hyperplane $L(d)$ is represented by the shaded region.}
\label{fig:3D}
\end{figure}

In Case III, $\ZZ=5k-6$ and $n^k_v - \nf= 3$.
Figure~\ref{fig:3D} illustrates Case III.  The three vertices must lie on the vertical edges of a triangular prism.  One vertex must be on a top-most vertex of this prism, and another vertex must be on a bottom-most vertex of the prism.
The final vertex can be located anywhere on the third vertical edge of the prism.  The algorithm uses 3 oracle calls to resolve the situation.  Figure \ref{fig:3D}, depicts one step of the inner loop within Case III.  The direction $d$ is selected such that $d^\top v =0$ at the three vertices represented by squared, and $d^\top(x^1+ u_k \e_k) = 1$ at the vertex represented by a circle.
The hyperplane $L(d)$ is represented by the shaded region.
Calling the oracle $\oracle(d)$ reveals the position of the vertex on the edge joining $x^1+ \ell_k \e_k$ and $x^1+ u_k \e_k$.  The final result is a total of $\ZZ+5 = 5k-1$ oracle calls.

\begin{figure}[htb!]
	
\begin{tikzpicture}[scale=0.32,
	node distance=5mm,
	cercle/.style={
                      rectangle,minimum size=3mm,rounded corners=3mm,
                      very thick,draw=black!50,
                      top color=white,bottom color=black!20,
                      font=\ttfamily}]
	\node at (0,9) [left]{\small $\R^{k-1}$};
	\node at (16.6,9) [right]{\small $\R^{k}$};
	\node at (0,7) [left]{\footnotesize $2k-3~$};
	\node at (0,4) [left]{\footnotesize $5k-8~$};
	\node at (0,1) [left]{\footnotesize $5k-6~$};
	\draw (14,7) -- (4,7);
	\draw (6,7) -- (14,4.5) ;
	\draw (14,4) -- (4,4);
	\draw (6,4) -- (14,1.5) ;
	\draw (14,1) -- (4,1) ;
	\node at (4,6.4) [right]{$^{+2}$};
	\node at (4,3.4) [right]{$^{+2}$};
	\node at (4,0.4) [right]{$^{+2}$};
	\node at (8,7.4) [right]{$^{+0}$};
	\node at (8,4.4) [right]{$^{+3}$};
	\node at (8,1.4) [right]{$^{+3}$};
	\node at (11,5.7) [rotate=-15]{$^{+0}$};
	\node at (11,2.7) [rotate=-15]{$^{+5}$};
	\node at (16.5,7) [right]{\footnotesize  $~(2k-3)+2=2k-1$};
	\node at (16.5,4) [right]{\footnotesize  $~\max\{(2k-3)+2, (5k-8)+2+3\} = 5k-3$};
	\node at (16.5,1) [right]{\footnotesize  $~\max\{(5k-8)+2+5, (5k-6)+2+3\} = 5k-1$};
	\node at (2,7)  [cercle]	{\scriptsize${n^{k-1}_v=1}$};
	\node at (2,4)  [cercle]	{\scriptsize${n^{k-1}_v=2}$};
	\node at (2,1)  [cercle]	{\scriptsize${n^{k-1}_v=3}$};
	\node at (15,7)  [cercle]	{\scriptsize${n^k_v=1}$};
	\node at (15,4)  [cercle]	{\scriptsize${n^k_v=2}$};
	\node at (15,1)  [cercle]	{\scriptsize${n^k_v=3}$};
\end{tikzpicture}

\caption{Number of oracle calls at iteration $k$ of Algorithm~\ref{Algo-Rn-nf3}.}
\label{fig-rec}
\end{figure}
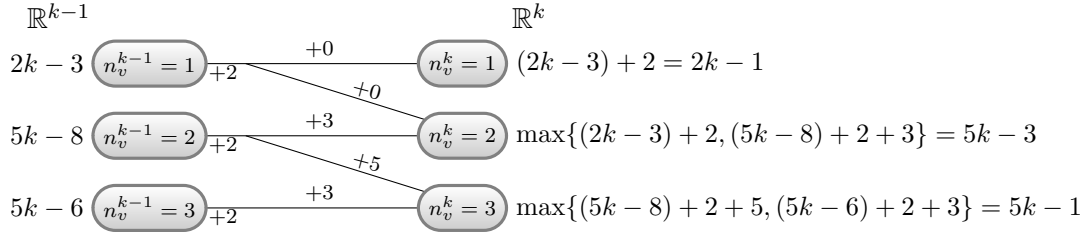

Figure~\ref{fig-rec} summarizes the total number of calls to the oracle to construct the vertices of the projection in $\R^k$ from those of $\R^{k-1}$. The left part of the figure lists $\ZZ$ the number of calls required to identify the $n_v^{k-1}$ vertices of the projection of $X$ onto $\R^{k-1}$.  Immediately after each of the three possible values of $n^{k-1}_v$ the symbol ``$+2$'' indicates that two oracle calls are made to compute $\ell_k$ and $u_k$.  The rest of the figure indicates the remaining number of oracle calls.

For example, if $n^k_v =3$ then  $n^{k-1}_v$ is necessarily equal to $2$ or $3$.
If $n^{k-1}_n =2$ then $5$ more calls are necessary (by Case II),  and if $n^{k-1}_n =3$ then $3$ more calls are necessary (by Case III). The maximal value between $(5k-8)+2 + 4$ and $(5k-6)+2+3$ constitutes an upper bound on the number of oracle calls  when $n^k_v=3$.\qed

\section{Discussion}\label{Conclusions}

We have studied the question of how to reconstruct a polyhedral subdifferential using directional derivatives.
By reformulating the question as the reconstruction of an arbitrary polyhedral $X$ set using an oracle $\oracle : \R^n \mapsto \R$ that returns $\displaystyle \oracle(d) = \max_{v \in X} v^\top d$, we observed that the question is closely linked to Geometric Probing.

We have developed a number of algorithms that provide methods to reconstruct a polyhedral subdifferential using directional derivatives in various situations.  However, many situations remain as open questions.  Table \ref{tab-summary} summarizes the results in this paper.

\begin{table}[ht]
\begin{center}
\begin{tabular}{|c||l|c|l|}
	\hline
	Space &Nb vertices and bound & Nb calls &  Source \\ \hline \hline
	\multirow{2}{1.5cm}{$~\R^1$}
		&	$\nf=\nv=1$	&	1 &	Section \ref{OneDim} \\
		&	$\nf \geq 2$	&	2 &	Section \ref{OneDim} \\
	\hline
	\multirow{5}{1.5cm}{$~\R^2$}
		&	$\nf=\nv=1$		&	$2$ 	&	Subsection \ref{nDim-nf1} \\
		&	$\nf>\nv=1$	&	$3$ 	&	Algorithm \ref{Algo-exact} with Theorem~\ref{Th:exactbounds} \\
		&	$\nf=\nv=2$		&	$5$ 	&	Algorithm \ref{Algo-exact} with Corollary~\ref{Cor:nf2} \\
		&	$\nf=\nv$			&	$3\nv$ &	Algorithm \ref{Algo-exact} with Theorem~\ref{Th:exactbounds}  \\
		&	$\nf>\nv$			&	$3\nv+1$ &	Algorithm \ref{Algo-exact} with Theorem~\ref{Th:exactbounds}  \\
	\hline
	\multirow{6}{1.5cm}{$\begin{array}{l}\R^n \end{array}$}
		&	$\nf=\nv=1$		&	$n$ 	&	Subsection \ref{nDim-nf1}  \\
		&	$\nf=2$, $\nv=1$  	&	$2n$ 	&	Algorithm \ref{Algo-Rn-nf2} with Theorem~\ref{Th-Rn-nf2} \\
		&	$\nf=2$, $\nv=2$		&	$3n-1$ 	&	Algorithm \ref{Algo-Rn-nf2} with Theorem~\ref{Th-Rn-nf2}  \\
		&	$\nf=3$, $\nv=1$	&	$2n-1$ 	&	Algorithm \ref{Algo-Rn-nf3} with Theorem~\ref{Th:Rn-nf3}  \\
		&	$\nf=3$, $\nv=2$	&	$5n-3$ 	&	Algorithm \ref{Algo-Rn-nf3} with Theorem~\ref{Th:Rn-nf3}  \\
		&	$\nf=3$, $\nv=3$		&	$5n-1$ 	&	Algorithm \ref{Algo-Rn-nf3} with Theorem~\ref{Th:Rn-nf3}  \\
	\hline
\end{tabular}
\end{center}
\caption{Maximum number of calls to the oracle to identify $X_v \subset \R^n$.}\label{tab-summary}
\end{table}%

From Table \ref{tab-summary} we see that the problem can be considered solved in $\R^1$ and $\R^2$.  However, in $\R^n$ for an arbitrary bound on the number of vertices $\nf$ the problem is still open.

This research is inspired in part by recent techniques that create approximate subdifferentials by using a collection of
approximate directional derivatives \cite{bagirov-2003, bagirov-karasozen-sezer-2008, Kiwiel-2010}.  As such, a natural research direction in this field is to examine how to adapt these algorithms to inexact oracles.   That is an oracle $\oracle^\varepsilon : \R^n \mapsto \R$ that returns $\displaystyle \oracle^\varepsilon(d) = \max_{v \in X} v^\top d + \xi$, where $\xi$ is an unknown error term bounded by $|\xi| < \varepsilon$. Some of the algorithms in this paper trivially adapt to this setting.
Specifically, the algorithm of Section \ref{OneDim} (for $\R^1$) and the algorithm of Subsection \ref{nDim-nf1} (for $\nf=1$) also work for inexact oracles and the error bounds are trivial to calculate.  However, the more interesting algorithms (Algorithm \ref{Algo-exact}, \ref{Algo-Rn-nf2} and \ref{Algo-Rn-nf3}) are not so trivial to adapt.

\begin{figure}[ht]
\begin{center}\includegraphics[width=8.7cm]{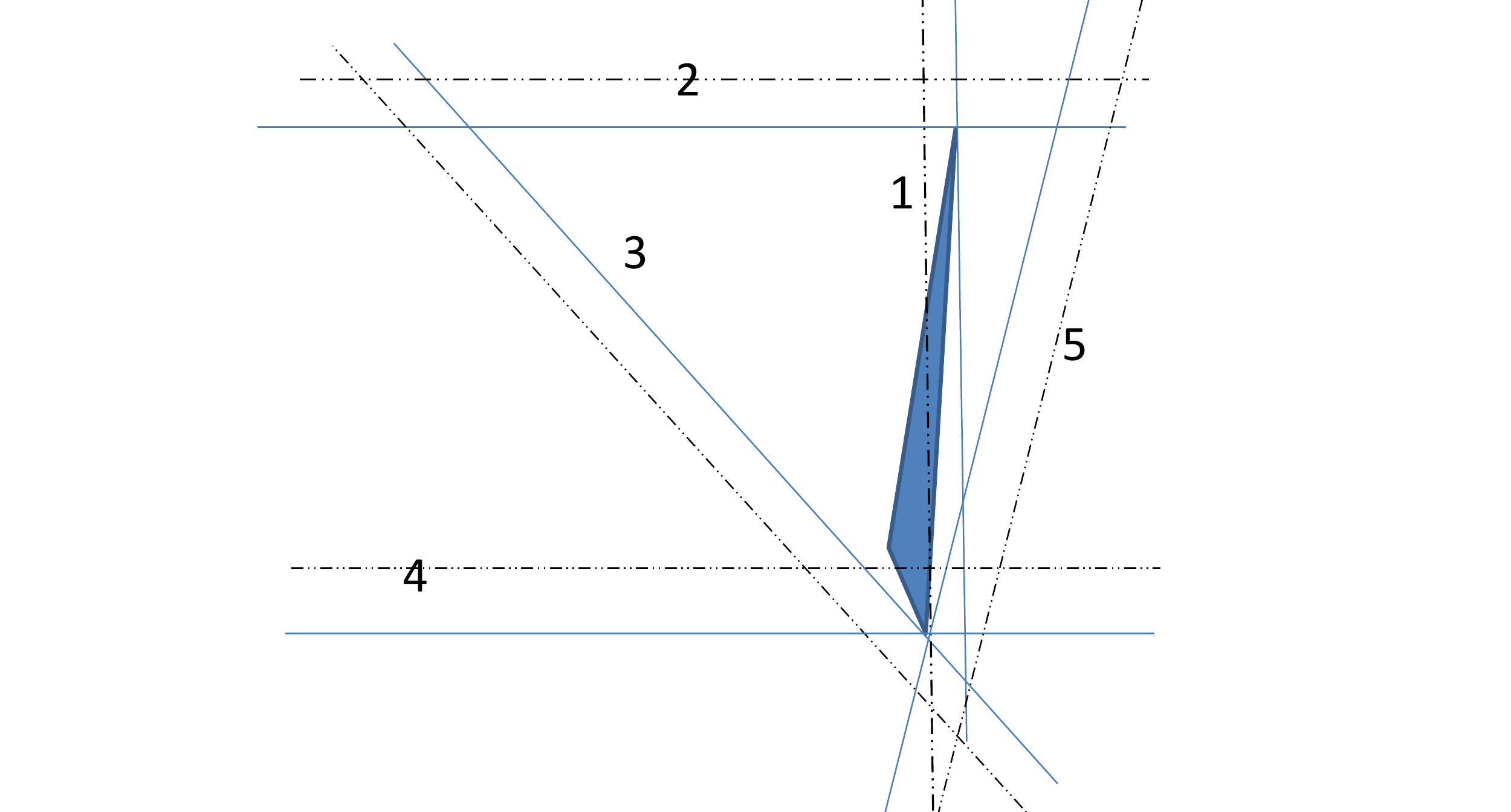}\end{center}
\caption{Possible outcome of Algorithm \ref{Algo-exact} when run with an inexact oracle.  After 5 oracle calls the algorithm stalls, as no new vertices are discovered.}\label{fig:inexact}
\end{figure}

We conclude this paper with Figure \ref{fig:inexact}, which demonstrates a potential problematic outcome of Algorithm \ref{Algo-exact} if run using an inexact oracle as if it were exact.  The continuous lines represent the hyperplane generated by an exact oracle, and the dotted ones are generated by an inexact one.  In this example, the fifth oracle call is incompatible with  Algorithm \ref{Algo-exact}, as no new vertices are discovered.

\subsection*{Acknowledgements} The authors are appreciative to the two anonymous referees for their time reviewing and constructive comments pertaining to this paper. 

\bibliographystyle{plain}
\bibliography{Bibliography}

\begin{thebibliography}{10}

\bibitem{AuDe03a}
C.~Audet and J.E. {Dennis, Jr.}
\newblock Analysis of generalized pattern searches.
\newblock {\em SIAM Journal on Optimization}, 13(3):889--903, 2003.

\bibitem{AuDe2006}
C.~Audet and J.E. {Dennis, Jr.}
\newblock Mesh adaptive direct search algorithms for constrained optimization.
\newblock {\em SIAM Journal on Optimization}, 17(1):188--217, 2006.

\bibitem{bagirov-2003}
A.M. Bagirov.
\newblock Continuous subdifferential approximations and their applications.
\newblock {\em Journal of Mathematical Sciences (New York)}, 115(5):2567--2609,
  2003.
\newblock Optimization and related topics, 2.

\bibitem{bagirov-karasozen-sezer-2008}
A.M. Bagirov, B.~Karas{\"o}zen, and M.~Sezer.
\newblock Discrete gradient method: derivative-free method for nonsmooth
  optimization.
\newblock {\em Journal of Optimization Theory and Applications},
  137(2):317--334, 2008.

\bibitem{Bernstein1986}
H.J. Bernstein.
\newblock Determining the shape of a convex n-sided polygon by using 2n + k
  tactile probes.
\newblock {\em Information Processing Letters}, 22(5):255 -- 260, 1986.

\bibitem{ColeYap84}
R.~Cole and C.K. Yap.
\newblock Shape from probing.
\newblock {\em Journal of Algorithms}, 8(1):19 -- 38, 1987.

\bibitem{DeRu1995}
V.F. Demyanov and A.M. Rubinov.
\newblock {\em Constructive nonsmooth analysis}.
\newblock Approximation and Optimization, vol. 7. Verlag Peter Lang, Frankfurt,
  1995.

\bibitem{Dobkin86}
D.~Dobkin, H.~Edelsbrunner, and C.K. Yap.
\newblock Probing convex polytopes.
\newblock In {\em Proceedings of the Eighteenth Annual ACM Symposium on Theory
  of Computing}, STOC '86, pages 424--432, New York, NY, USA, 1986. ACM.

\bibitem{FiKe2004a}
D.E. Finkel and C.T. Kelley.
\newblock Convergence analysis of the {DIRECT} algorithm.
\newblock Technical Report CRSC-TR04-28, Center for Research in Scientific
  Computation, 2004.

\bibitem{FiKe09}
D.E. Finkel and C.T. Kelley.
\newblock Convergence analysis of sampling methods for perturbed lipschitz
  functions.
\newblock {\em Pacific Journal of Optimization}, 5(2):339--350, 2009.

\bibitem{Hare-Nutini-2013}
W.~Hare and {J. Nutini}.
\newblock A derivative-free approximate gradient sampling algorithm for finite
  minimax problems.
\newblock {\em Computational Optimization and Applications}, 56(1):1--38, 2013.

\bibitem{hare-macklem-2012}
W.~Hare and M.~Macklem.
\newblock Derivative-free optimization methods for finite minimax problems.
\newblock {\em Optimization Methods and Software}, iFirst:1--13, 2011.

\bibitem{Hiriart-Urruty-Lemarechal-1993b}
J.B. Hiriart-Urruty and C.~Lemar{\'e}chal.
\newblock {\em Convex Analysis and Minimization Algorithms. {II}}, volume 306
  of {\em Grundlehren der Mathematischen Wissenschaften [Fundamental Principles
  of Mathematical Sciences]}.
\newblock Springer-Verlag, Berlin, 1993.
\newblock Advanced theory and bundle methods.

\bibitem{Imiya2012}
A.~Imiya and K.~Sato.
\newblock Shape from silhouettes in discrete space.
\newblock In E.V. Brimkov and P.R. Barneva, editors, {\em Digital Geometry
  Algorithms: Theoretical Foundations and Applications to Computational
  Imaging}, pages 323--358. Springer Netherlands, Dordrecht, 2012.

\bibitem{Kiwiel-2010}
K.C. Kiwiel.
\newblock A nonderivative version of the gradient sampling algorithm for
  nonsmooth nonconvex optimization.
\newblock {\em SIAM Journal on Optimization}, 20(4):1983--1994, 2010.

\bibitem{Lindenbaum92}
M.~Lindenbaum and A.~Bruckstein.
\newblock Parallel strategies for geometric probing.
\newblock {\em Journal of Algorithms}, 13(2):320 -- 349, 1992.

\bibitem{rockafellar-wets-1998}
R.T. Rockafellar and R.J.-B. Wets.
\newblock {\em Variational Analysis}, volume 317 of {\em Grundlehren der
  Mathematischen Wissenschaften [Fundamental Principles of Mathematical
  Sciences]}.
\newblock Springer-Verlag, Berlin, 1998.

\bibitem{Romanik95}
K.~Romanik.
\newblock Geometric probing and testing - a survey.
\newblock Technical Report 95-42, DIMACS Technical Report, 1995.

\bibitem{Li1988}
R.L. Shuo-Yen.
\newblock Reconstruction of polygons from projections.
\newblock {\em Information Processing Letters}, 28(5):235 -- 240, 1988.

\bibitem{Skiena89}
S.S. Skiena.
\newblock Problems in geometric probing.
\newblock {\em Algorithmica}, 4(1):599--605, 1989.

\bibitem{Tiwary2008}
H.R. Tiwary.
\newblock {\em Complexity of some polyhedral enumeration problems}.
\newblock PhD thesis, Universit\"at des Saarlandes, Germany, 2008.

\end{thebibliography}

\end{document}